# Elastically-isotropic open-cell minimal surface shell lattices with superior stiffness via variable thickness design

Qingping Ma [a], Lei Zhang [a], Junhao Ding [b], Shuo Qu [b], Jin Fu [c], Mingdong Zhou [d], Ming Wang Fu [c], Xu Song [b], Michael Yu Wang [a, *]

[a] *Department of Mechanical and Aerospace Engineering, Hong Kong University of Science and Technology, Kowloon, Hong Kong, China*

[b] *Department of Mechanical and Automation Engineering, Chinese University of Hong Kong, Shatin, Hong Kong, China*

[c] *Department of Mechanical Engineering, Hong Kong Polytechnic University, Kowloon, Hong Kong, China*

[d] *Shanghai Key Laboratory of Digital Manufacture for Thin-walled Structures, School of Mechanical Engineering, Shanghai Jiao Tong University, Shanghai, China*

[*] Corresponding author: Michael Yu Wang (mywang@ust.hk)

# Abstract

Triply periodic minimal surface (TPMS) shell lattices are attracting increasingly attention due to their unique combination of geometric and mechanical properties, and their open-cell topology. However, uniform thickness TPMS shell lattices are usually anisotropic in stiffness, namely having different Young's moduli along different lattice directions. To reduce the elastic anisotropy, we propose a family of variable thickness TPMS shell lattices with isotropic stiffness designed by a strain energy-based optimization algorithm. The optimization results show that all the six selected types of TPMS lattices can be made to achieve isotropic stiffness by varying the shell thickness, among which N14 can maintain over 90% of the Hashin-Shtrikman upper bound of bulk modulus. All the optimized shell lattices exhibit superior stiffness properties and significantly outperform elastically-isotropic truss lattices of equal relative densities. Both uniform and optimized types of N14 shell lattices along [100], [110] and [111] directions are fabricated by the micro laser powder bed fusion techniques with stainless steel 316L and tested under quasi-static compression loads. Experimental results show that the elastic anisotropy of the optimized N14 lattices is reduced compared to that of the uniform ones. Large deformation compression results reveal different failure deformation behaviors along different directions. The [100] direction shows a layer-by-layer plastic buckling failure mode, while the failures along [110] and [111] directions are related to the shear deformation. The optimized N14 lattices possess a reduced anisotropy of plateau stresses and can even attain nearly isotropic energy absorption capacity.

**Keywords:** Elastic isotropy; Open-cell shell lattices; Triply periodic minimal surface; Laser powder bed fusion; Mechanical properties

# 1. Introduction

Mechanical metamaterials are artificial structures whose mechanical properties are determined by the downscale micro-architectures [1]. By tailoring the architectures, mechanical metamaterials can obtain unconventional mechanical properties unachievable by bulk materials, such as negative Poisson's ratio (auxetics) [2], negative compressibility transitions [3], pentamode properties [4], and maintain high stiffness and strength with low relative densities [1, 5]. Recent advances in additive manufacturing (AM) technologies facilitate the realization of mechanical metamaterials with complex downscale micro-architectures [6, 7].

In recent decades, numerous research efforts have been made to improve the mechanical performance of metamaterials in terms of architecture design. The stochastic foams first appeared as a class of artificial porous elastically-isotropic metamaterials, mainly for energy absorption applications [1, 8]. The subsequent emerging open-cell truss lattices can outperform stochastic foams of equal relative densities in stiffness and strength via the stretching-dominated design concept, and therefore receive considerable attention [9-11]. However, even the optimal stretching-dominated elastically-isotropic truss lattices can only reach nearly 1/3 of the Hashin-Shtrikman (HS) upper bounds in low relative density limits [12-14]. Recent studies have found a family of stretching-dominated elastically-isotropic plate lattices whose stiffness properties can reach the HS upper bounds in low relative density limits [15-17]. However, their closed-cell topology complicates the manufacturing and therefore additional holes need to be added for resin or powder removal in the AM process [16-18]. In addition, many applications such as heat exchangers [19], biomedical implants [20, 21] and supercapacitors [22-24] require smooth open channels for mass and heat transfer. Thus, it remains a significant research task to find open-cell metamaterials with superior mechanical properties for multifunctional applications.

Other than stochastic foams, truss lattices and plate lattices, shell lattices provide another opportunity of mechanical metamaterials. They are now attracting increasingly attention for their smooth, periodic and non-intersecting structures to avoid stress concentration and their open-cell topology. The majority of studies on shell lattices are based on the well-known triply periodic minimal surface (TPMS) [25]. The stiffness and strength of TPMS shell lattices are proven to exceed those of truss lattices in a wide range of relative densities [21, 26]. Especially, the bulk moduli of TPMS shell lattices are numerically confirmed to approach the HS upper bound in low relative density limits [27]. Functionally graded designs of TPMS shell lattices are further developed for enhanced energy absorption capacity [28], exceptional fatigue resistance [29], bone-mimicking transport properties [30], and are widely adopted for tissue engineering scaffolds. Yoo [31] proposed a heterogeneous modeling method to design TPMS scaffolds with controllable porosity, while Vijayavenkataraman et al. [32] further developed a parametric optimization approach for sheet TPMS scaffolds to satisfy multifunctional requirements including stiffness, porosity and pore size.

The studies on shell lattices above are related to TPMS shell lattices with uniform thickness, which are usually elastically-anisotropic – possessing different stiffness properties along different lattice directions. The stiffness properties along some directions may be much below the HS upper bounds [26, 33], while a deviation of the loading direction may lead to a large deformation or even catastrophic failure of the anisotropic lattices, which limits their applications once the primary loading direction is unknown. Therefore, the elastically-isotropic shell lattices are preferred for applications in such cases, while several attempts have so far been made to achieve macroscopic elastic isotropy. Bonatti et al. [26, 33] proposed a biased functional based shape optimization for generating new families of elastically-isotropic open-cell shell lattices. Soyarslan et al. [34] developed a parametric shape control method based on the periodic nodal surface to achieve elastically-isotropic shell lattices. Chen et al. [27] and Deng et al. [35] proposed several elastically-isotropic hybrid lattices by combinations of different types of TPMS shell lattices or combinations of plate lattices and TPMS shell lattices, while such combinations change the open-cell topology to closed-cell and additional holes have to be added to make the lattices manufacturable. Callens et al. [36] proposed a parametric metamaterial design strategy and achieved elastically-isotropic TPMS shell lattices through multi-material hyperbolic tiling. The related work either changes the shape of the shell mid-surface away from minimal surface [26, 33, 34], or divides the lattice into several components with different lattice types or material properties [27, 35, 36]. Given the existing simulation results that shell lattices with minimal surface geometry can approach the HS upper bound of bulk modulus in low relative density limits [27], it turns out to be an inspiring goal to find open-cell, single-material and elastically-isotropic shell lattices which maintain the minimal surface as the shell mid-surface to achieve superior stiffness.

In this work, a family of elastically-isotropic open-cell variable thickness TPMS shell lattices are designed. The elastic anisotropy of a TPMS shell lattice is tuned by smoothly varying its shell thickness with a strain energy-based numerical homogenization and optimization procedure. As a result, six types of elastically-isotropic TPMS shell lattices are generated, among which N14 is found to maintain over 90% of the HS upper bound of bulk modulus. All the optimized shell lattices exhibit superior stiffness properties and significantly outperform elastically-isotropic truss lattices of equal relative densities. The elastically-isotropic N14 shell lattices and the conventional counterparts of uniform thickness are fabricated with the micro laser powder bed fusion (LPBF) techniques in stainless steel 316L (SS316L). Quasi-static compression tests are conducted to investigate the Young's moduli, strength, and failure deformation behaviors along the three principal lattice directions ([100], [110] and [111]) and to validate the elastic isotropy of the variable thickness designs.

## 2. Elastically-isotropic design of TPMS shell lattices

### 2.1 TPMS shell lattices

In mathematics, a TPMS is a triply periodic surface with constant zero mean curvature. The Enneper-Weierstrass representation gives the accurate coordinates of the TPMS [37], while the Weierstrass function is only known for several types of TPMS. Alternatively, the TPMS is usually approximated by the periodic nodal surface with enough Fourier terms, which enables much easier generation of TPMS but the zero mean curvature may not be guaranteed [38]. Surface Evolver is a powerful tool for generating discrete TPMS [39, 40], minimizing the surface tension energy of a surface with the prescribed boundary conditions and eventually leading to a TPMS with good accuracy. In this work, an accurate TPMS is generated in stereolithography (STL) format with the use of Surface Evolver and then used as the shell mid-surface for finite element analysis (FEA) and optimization to achieve elastically-isotropic shell lattices via variable thickness design.

This work focuses on the TPMS with cubic symmetry, i.e., with reflectional symmetry about the three mid-planes and orthogonal rotational symmetry about the three axes. The former symmetry enables a unit cell to be decomposed into 8 equal parts, while the latter symmetry further divides each 1/8 cell into 6 smaller equal parts. As a result, the 1/48 domain (quadrirectangular tetrahedron) is a fundamental unit of the unit cell, based on which the 1/8 cell and unit cell can be generated by mirror operations. Six types of TPMS are shown in Fig. 1, including the N14, OCTO, I-graph-wrapped package (IWP), Face-centered cubic rhombic dodecahedron (FRD), Neovius (N) and Primitive (P) surfaces, to illustrate the 1/48 fundamental unit (red), 1/8 cell (blue) and unit cell (grey), respectively. The 1/8 cell can be generated by mirroring the 1/48 fundamental unit about the neighboring edges, while the unit cell can be generated by mirroring the 1/8 cell about the three mid-planes.

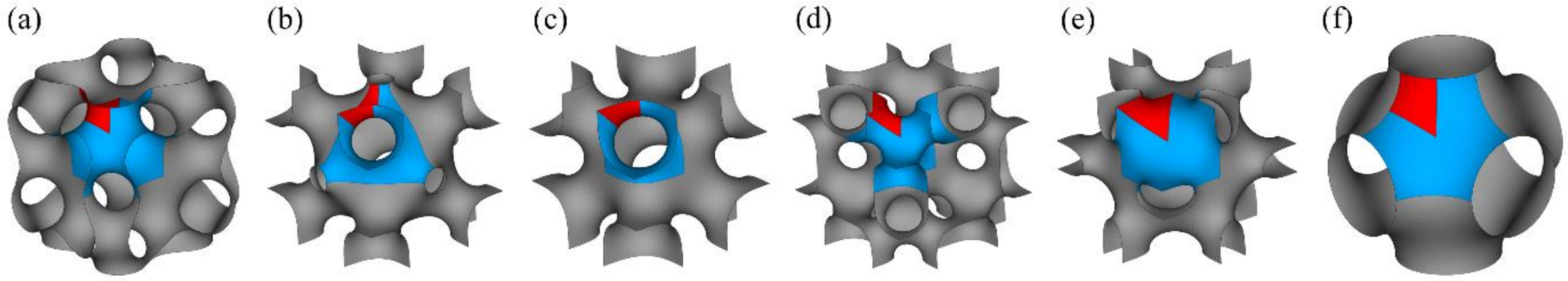

Fig. 1 Six types of TPMS, including: (a) N14, (b) OCTO, (c) IWP, (d) FRD, (e) N, (f) P, with red, blue and grey colors representing the 1/48 fundamental unit, 1/8 cell and unit cell, respectively.

## 2.2 Elasticity of cubic symmetric lattices

The macroscopic constitutive relation of cubic symmetric lattices is represented as [41]:

$$\boldsymbol{\sigma} = \boldsymbol{C}\boldsymbol{\varepsilon}, \tag{1}$$

where the macroscopic stress $\boldsymbol{\sigma}$, homogenized elasticity tensor $\boldsymbol{C}$ and macroscopic strain $\boldsymbol{\varepsilon}$ are expressed in the Voigt notation as:

$$\boldsymbol{\sigma} = \begin{pmatrix} \sigma_{11} \\ \sigma_{22} \\ \sigma_{33} \\ \sigma_{23} \\ \sigma_{13} \\ \sigma_{12} \end{pmatrix}, \quad \boldsymbol{C} = \begin{pmatrix} c_{11} & c_{12} & c_{12} & 0 & 0 & 0 \\ c_{12} & c_{11} & c_{12} & 0 & 0 & 0 \\ c_{12} & c_{12} & c_{11} & 0 & 0 & 0 \\ 0 & 0 & 0 & c_{44} & 0 & 0 \\ 0 & 0 & 0 & 0 & c_{44} & 0 \\ 0 & 0 & 0 & 0 & 0 & c_{44} \end{pmatrix}, \quad \boldsymbol{\varepsilon} = \begin{pmatrix} \varepsilon_{11} \\ \varepsilon_{22} \\ \varepsilon_{33} \\ \gamma_{23} \\ \gamma_{13} \\ \gamma_{12} \end{pmatrix}, \tag{2}$$

where $c_{11}$, $c_{12}$ and $c_{44}$ are three independent elasticity constants which completely define the stiffness properties of cubic symmetric lattices. The elastic anisotropy is typically assessed with the Zener anisotropic index $\xi$:

$$\xi = \frac{2c_{44}}{c_{11} - c_{12}}, \tag{3}$$

where the Zener anisotropic index $\xi = 1$ corresponds to elastic isotropy.

Typically, three load cases are adopted for evaluation of the three elasticity constants, including the uniaxial strain loading, hydrostatic loading and pure shear loading. Since the value of the imposed strain has no influence on the stiffness properties of the lattices in the linear elastic model, a unit value is adopted for simplicity, i.e., $\boldsymbol{\varepsilon}^{(U)} = [1,0,0,0,0,0]^T$, $\boldsymbol{\varepsilon}^{(T)} = [1,1,1,0,0,0]^T$ and $\boldsymbol{\varepsilon}^{(S)} = [0,0,0,0,0,1]^T$, where the superscripts (*U*), (*T*) and (*S*) correspond to the uniaxial strain loading, hydrostatic loading and pure shear loading, respectively. Based on Eq. (1) and (2), the three elasticity constants are obtained as:

$$\begin{aligned} &c_{11} = \sigma_{11}^{(U)}, \\ &c_{12} = \sigma_{22}^{(U)} = \sigma_{33}^{(U)}, \\ &c_{44} = \sigma_{12}^{(S)}, \\ &c_{11} + 2c_{12} = \sigma_{11}^{(T)} = \sigma_{22}^{(T)} = \sigma_{33}^{(T)}, \end{aligned} \tag{4}$$

where $\sigma_{11}^{(U)}$, $\sigma_{22}^{(U)}$, $\sigma_{33}^{(U)}$, $\sigma_{12}^{(S)}$, $\sigma_{11}^{(T)}$, $\sigma_{22}^{(T)}$ and $\sigma_{33}^{(T)}$ are the macroscopic stress components of the corresponding load cases. Although only two load cases (uniaxial strain loading and pure shear loading) are enough for determining all the three elasticity constants, the hydrostatic loading is necessary for the sensitivity analysis of the optimization problem to achieve elastic isotropy, as detailed in Section 2.4. In addition, the hydrostatic loading serves as a verification of the two elasticity constants $c_{11}$ and $c_{12}$ obtained from the uniaxial strain loading. To evaluate the macroscopic stresses $\sigma_{11}$, $\sigma_{22}$, $\sigma_{33}$ and $\sigma_{12}$ in Eq. (4), two expressions are adopted based on Hill's principle [41]:

$$\boldsymbol{\sigma} = \frac{1}{A}\int_A d\boldsymbol{F}, \tag{5a}$$

$$\boldsymbol{\sigma} = \frac{1}{V}\int_V \boldsymbol{\sigma}_\mu dV, \tag{5b}$$

where $A$, $V$ are the end face area and volume of the selected representative volume element, $\boldsymbol{F}$ is the reaction force on the end faces, while $\boldsymbol{\sigma}_\mu$ is the microscopic stress. Eq. (5a) represents the average stress over the end faces and is called the macroscopic average stress, while Eq. (5b) represents the average stress over the entire volume of the representative element and is called the microscopic average stress. The Hill's principle reveals that the macroscopic average stress is equal to the microscopic average stress [41].

As shown in Section 2.5, a strain energy-based optimization algorithm is proposed in this work for achieving elastically-isotropic shell lattices, while the elasticity constants in Eq. (4) are in terms of the macroscopic stresses. For derivation of the sensitivity, the elasticity constants are rewritten in terms of the strain energy. To achieve this, the strain energies of the three load cases are first obtained in terms of the macroscopic stresses as:

$$e^{(U)} = \frac{1}{2}\sigma_{11}^{(U)}V, \tag{6a}$$

$$e^{(T)} = \frac{1}{2}[\sigma_{11}^{(T)} + \sigma_{22}^{(T)} + \sigma_{33}^{(T)}]V, \tag{6b}$$

$$e^{(S)} = \frac{1}{2}\sigma_{12}^{(S)}V, \tag{6c}$$

Based on Eq. (4), the strain energies in Eq. (6a) ~ (6c) are rewritten in terms of the elasticity constants as:

$$e^{(U)} = \frac{1}{2}c_{11}V, \tag{7a}$$

$$e^{(T)} = \frac{3}{2}(c_{11} + 2c_{12})V, \tag{7b}$$

$$e^{(S)} = \frac{1}{2}c_{44}V, \tag{7c}$$

Thus, the elasticity constants in Eq. (4) are rewritten in terms of the strain energy as:

$$c_{11} = \frac{2}{V}e^{(U)}, \tag{8a}$$

$$c_{12} = \frac{1}{V}[\frac{1}{3}e^{(T)} - e^{(U)}], \tag{8b}$$

$$c_{44} = \frac{2}{V}e^{(S)}, \tag{8c}$$

Based on the three elasticity constants, the stiffness properties of cubic symmetric lattices are obtained as [27, 34]:

$$E = \frac{(c_{11} - c_{12})(c_{11} + 2c_{12})}{c_{11} + c_{12}}, \tag{9a}$$

$$G = c_{44}, \tag{9b}$$

$$K = \frac{c_{11} + 2c_{12}}{3}, \tag{9c}$$

$$\nu = \frac{c_{12}}{c_{11} + c_{12}}, \tag{9d}$$

where $E$, $G$, $K$ and $\nu$ denote the Young's, shear and bulk moduli and Poisson's ratio of the lattices, respectively. For elastically-isotropic lattices with $\xi = 1$, the HS upper bounds define the maximum stiffness properties [15-17, 42]:

$$\frac{E}{E_s} = \frac{2\bar{\rho}(5\nu_s - 7)}{13\bar{\rho} + 12\nu_s - 2\bar{\rho}\nu_s - 15\bar{\rho}\nu_s^2 + 15\nu_s^2 - 27}, \tag{10a}$$

$$\frac{K}{K_s} = \frac{4\bar{\rho}G_s}{4G_s + 3K_s(1 - \bar{\rho})}, \tag{10b}$$

where $E_s$, $G_s$, $K_s$ and $\nu_s$ denote the Young's, shear and bulk moduli and Poisson's ratio of the constituent material, respectively, while $\bar{\rho}$ denotes the relative density of the lattices.

## 2.3 Initial yield properties

Furthermore, the effects of thickness variations on the initial yield properties of the lattices are studied, in which linear elastic FEA under three load cases are adopted for comparisons, including the uniaxial stress loading, hydrostatic loading and pure shear loading. Given the macroscopic constitutive relation in Eq. (1) and (2), the following macroscopic strain is adopted [14, 16]:

$$\boldsymbol{\varepsilon}^{(US)} = \left[1, -\frac{c_{12}}{c_{11} + c_{12}}, -\frac{c_{12}}{c_{11} + c_{12}}, 0, 0, 0\right]^T, \tag{11}$$

to induce a macroscopic uniaxial stress state of:

$$\boldsymbol{\sigma}^{(US)} = \left[\frac{(c_{11} - c_{12})(c_{11} + 2c_{12})}{c_{11} + c_{12}}, 0, 0, 0, 0, 0\right]^T, \tag{12}$$

where the superscript ($US$) represents the uniaxial stress loading, while the other two load cases are identical to those discussed above. The macroscopic stresses $\Sigma$ for the three load cases are defined as:

$$\Sigma^{(US)} = \sigma_{11}^{(US)}, \quad \Sigma^{(T)} = \frac{1}{3}[\sigma_{11}^{(T)} + \sigma_{22}^{(T)} + \sigma_{33}^{(T)}], \quad \Sigma^{(S)} = \sigma_{12}^{(S)}, \tag{13}$$

based on which the normalized initial yield stresses of the lattices are derived as [14]:

$$(\frac{\Sigma_s^y}{\rho\sigma_s^y})^{(i)} = \frac{\Sigma^{(i)}}{\rho(\sigma_{vM}^{\max})^{(i)}}, \quad (i = US, T, S), \tag{14}$$

where $\Sigma_s^y$ denotes the initial yield stress of the lattices, $\sigma_s^y$ denotes the yield strength of the constituent material, while $\sigma_{vM}^{\max}$ denotes the maximum von Mises stress of the lattices under the specified load case.

### 2.4 FEA implementation

To evaluate the three elasticity constants, FEA under the three corresponding load cases are performed. A representative unit cell is typically adopted for FEA, with the periodic boundary conditions imposed on the node pairs of opposite end faces [34, 41, 43-45]:

$$\begin{aligned} &u_i^+ - u_i^- = \varepsilon_{ij}(X_j^+ - X_j^-), \quad (i = 1, 2, 3), \\ &\theta_i^+ - \theta_i^- = 0, \quad (i = 1, 2, 3), \end{aligned} \tag{15}$$

where the superscripts $+$ and $-$ represent the corresponding nodes on the opposite end faces, $j$ is a dummy index of the Einstein summation convention, $\boldsymbol{u}$ and $\boldsymbol{\theta}$ denote the translational and rotational displacement vectors, respectively, while $\boldsymbol{X}$ denotes the locations of the nodes.

Due to the reflectional symmetry of the TPMS considered in this work, the 1/8 cell is adopted for FEA to save the computational cost without loss of accuracy [41, 46]. The P-surface is adopted to illustrate the boundary conditions in Fig. 2, in which the blue and grey colors represent the surface within the 1/8 cell and unit cell, respectively, while $x_0$, $x_-$ and $x_+$ represent the mid-plane, negative and positive end planes along $x$ direction, respectively. For the uniaxial strain loading and hydrostatic loading, the symmetry boundary conditions are imposed on the three mid-planes $x_0$, $y_0$ and $z_0$. For the pure shear loading, the anti-symmetry boundary conditions are imposed on the two mid-planes $x_0$ and $y_0$, while the symmetry boundary conditions are imposed on the other mid-plane $z_0$. The boundary conditions imposed on the three end planes $x_+$, $y_+$ and $z_+$ for all the three load cases can be derived by simplifying the expressions of the periodic boundary conditions in Eq. (15), as detailed in *Appendix A*.

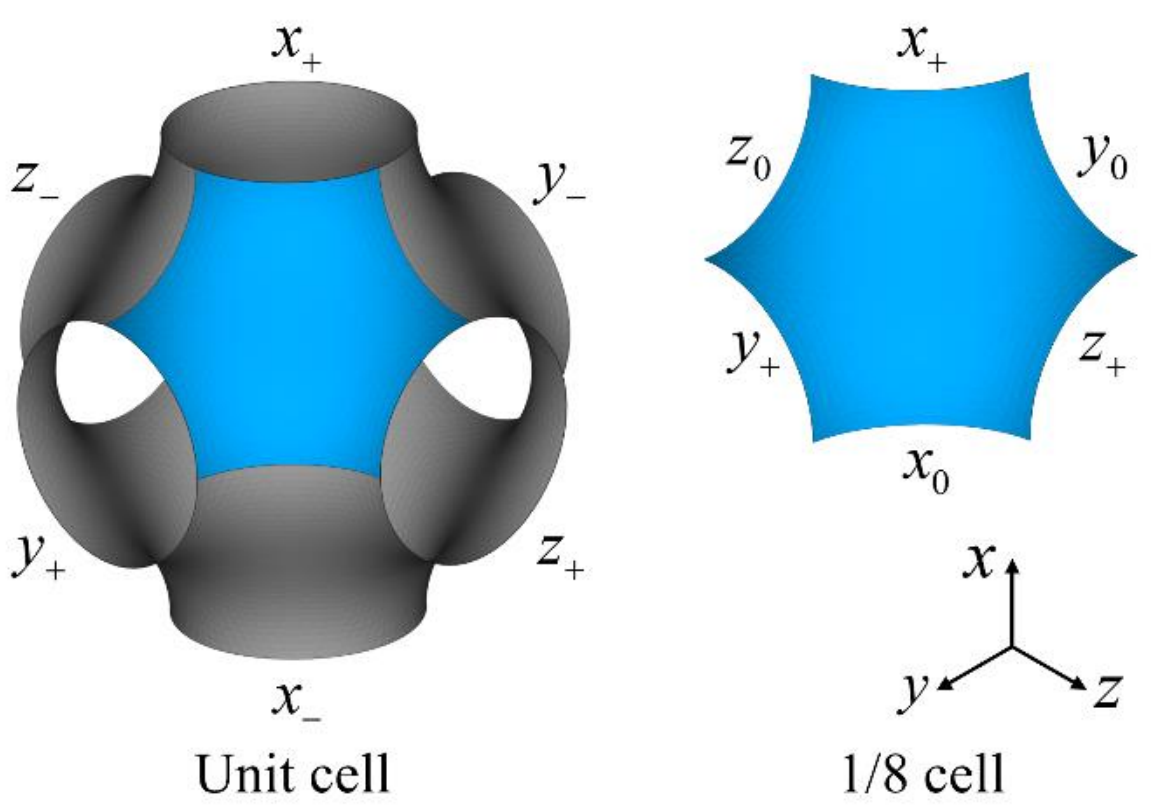

Fig. 2 Illustrations of the boundary conditions imposed on the unit cell and 1/8 cell by the TPMS P-surface.

This work adopts the ABAQUS STRI3 shell element for linear elastic analysis, based on the classical Kirchhoff-Love plate/shell theory. The constituent material is assumed to be isotropic in stiffness, with the Young's modulus $E_s = 190$ GPa and Poisson's ratio $\nu_s = 0.3$. For ease of FEA implementation, Eq. (5a) is adopted for evaluation of the macroscopic stress $\boldsymbol{\sigma}$, implemented as:

$$\boldsymbol{\sigma} = \frac{1}{A}\sum_{i=1}^{n} \boldsymbol{F}_i^R, \tag{16}$$

where $\boldsymbol{F}_i^R$ is the reaction force in the direction of the applied displacement, while $n$ is the number of nodes on the corresponding end faces. To achieve the macroscopic stiffness properties, the macroscopic stress is first evaluated by Eq. (16), based on which the strain energies of the three load cases are obtained by Eq. (6a) ~ (6c). The three elasticity constants are then evaluated by Eq. (8a) ~ (8c), based on which the stiffness properties of the lattices are achieved by Eq. (9a) ~ (9d), while the Zener anisotropic index is obtained by Eq. (3).

### 2.5 Elastically-isotropic TPMS shell lattices via variable thickness design

In this work, a family of elastically-isotropic TPMS shell lattices are proposed via variable thickness design, based on the following optimization problem:

$$\begin{cases} \min\limits_{\delta_i} J = \dfrac{1}{2}(\xi - 1)^2, \quad (i = 1, 2, \ldots, N), \\ s.t. \begin{cases} a(\boldsymbol{v}^{(j)}, \boldsymbol{w}^{(j)}) = l(\boldsymbol{w}^{(j)}), \quad \forall \boldsymbol{w}^{(j)} \in \boldsymbol{U}, \quad (j = U, T, S), \\ g(\delta) = \dfrac{1}{V}\displaystyle\int_M \delta dM - \rho^* = 0, \\ \delta_{\min} \le \delta_i \le \delta_{\max}, \quad (i = 1, 2, \ldots, N), \end{cases} \end{cases} \tag{17}$$

where the shell mid-surface is discretized into $N$ triangular elements, each of which is assigned with the thickness of $\delta_i$, while the upper and lower bounds of the thicknesses are $\delta_{\max}$ and $\delta_{\min}$, respectively. The objective function $J$ is defined in terms of the Zener anisotropic index $\xi$. The static equilibrium is given in the weak variational form in terms of the energy bilinear form $a(\boldsymbol{v}^{(j)},\boldsymbol{w}^{(j)})$ and load linear form $l(\boldsymbol{w}^{(j)})$, where $\boldsymbol{v}^{(j)}$ and $\boldsymbol{w}^{(j)}$ denote the displacement field and virtual displacement field of load case $j$, respectively, while $\boldsymbol{U}$ denotes the kinematically admissible displacement field. The constraint function $g(\delta)$ represents the volume constraint of the lattice, in which $M$ and $\rho^*$ denote the shell mid-surface and targeted relative density of the lattice, respectively.

This work focuses on the novel design of thin shell lattices with low relative densities for lightweight applications. The classical Kirchhoff-Love shell theory adopted in this work can accurately predict the mechanical properties of thin shells, while it may lead to inaccuracy for thick shells due to the neglection of transverse shear forces. As a result, the maximum thickness $\delta_{\max}$ is adopted as $(1/10)D$ for the thin shell assumption to hold, in which $D$ denotes the unit cell size, while the minimum thickness $\delta_{\min}$ should be appropriately selected for the optimization algorithm to yield elastically-isotropic designs which satisfy the manufacturing constraints of minimum feature size.

The sensitivity derivation and updating scheme of the shell thickness are detailed in *Appendix B*. The initial design is adopted as the uniform thickness TPMS shell lattices. In the optimization process, the 1/8 cell is taken as the FEA domain due to the reflectional symmetry, while the 1/48 fundamental unit is taken as the design domain since it is the smallest unit of simple cubic (SC) lattices. As discussed previously, the 1/8 cell consists of 6 fundamental domains, each of which can be generated as the mirror image of the neighboring fundamental domain. Since the strain energy distributions as well as the sensitivities of the 6 fundamental domains are not identical, an average operation of the sensitivity terms is performed in each iteration step to maintain cubic symmetry, i.e.,

$$\overline{\frac{\partial L}{\partial \delta_i}} = \frac{1}{6}\sum_{k=1}^{6}\left(\frac{\partial L}{\partial \delta_i}\right)^{(k)}, \quad (i = 1, 2, ..., N), \tag{18}$$

where the superscript ($k$) denotes the $k$-th fundamental domain within the 1/8 cell. The averaged sensitivity terms $\overline{\frac{\partial L}{\partial \delta_i}}$ are then adopted for updating the shell thickness by Eq. (B11) and (B12). The relative tolerances for convergence of the Zener anisotropic index and constraint function are set as 0.01, i.e., $|\xi - 1| \le 0.01$ and $|\frac{g(\delta)}{\rho^*}| \le 0.01$.

### 2.6 Variable-thickness offset operation

After the optimization iteration converges, a variable-thickness offset operation is performed by MATLAB programming to export the final design in STL format. First, the thickness of each node is taken as the average of the thicknesses of neighboring triangular elements. Then, the top and bottom surfaces are generated by offsetting every mid-surface node along its positive and negative normal directions for a half of the nodal thickness. A variable thickness shell lattice is then obtained by combing the top and bottom surfaces with the lateral surfaces.

## 3. Experimental methods

### 3.1 Micro-LPBF

To validate the stiffness properties, the uniform and variable thickness N14 shell lattices were fabricated on the *Hans' Laser M100* machine equipped with 500W IPG fiber laser (with the wavelength of 1.07 μm) with the beam diameter of 25 μm. The lens of the machine is custom-made with a standard SCANLAB galvanometer, especially the laser beam expander lens is custom-made to increase the incident laser diameter and to reduce the final laser spot size after the f-theta lens. All the lattices were fabricated using the gas atomized austenitic SS316L powder provided by Beijing AMC Powder Metallurgy Technology Co., Ltd. The particle size of the powder is within the range of 5 ~ 25 μm ($\text{D50} = 16.27\ \mu\text{m}$) and its chemical compositions are listed in Table 1. The laser power, scanning speed, hatch distance and layer thickness are 50 W, 1000 mm/s, 50 μm and 10 μm, respectively. A common scanning strategy was adopted to reduce the thermal stress by rotating the patterns between successive layers for a hatch angle of $67^\circ$. The process parameters have been optimized in our previous work [47] and the optical inspection on the cross sections of the bulk specimens revealed that the internal structures were dense and pore-free. Therefore, the lattice structures fabricated by the optimal parameters are in the full-densification condition. The as-printed lattices were removed from the baseplate by the wire electrical discharge machining and cleaned by the ethanol using the ultrasonic vibration. No heat treatment was conducted on the as-printed lattices, since the SS316L specimens are proven to possess high strength and ductility [48] due to the solidification-enabled cellular structures of dislocations [47, 49], while their mechanical properties will be reduced by the heat treatment process [50].

The micro-LPBF fabricated SS316L has the density of $\rho_s = 8.03\ \text{g/cm}^3$ and Young's modulus of $E_s = 190\ \text{GPa}$, based on the in-house standard tensile tests of tensile specimens, in which the loading direction is within the platform ($x$-$y$) plane.

Table 1 Chemical compositions of SS316L

| Element | Cr | Ni | Mo | Mn | Si | P | S | C | O | Fe |
|---|---|---|---|---|---|---|---|---|---|---|
| wt. % | 16.88 | 13.6 | 2.7 | 0.54 | 0.41 | 0.012 | 0.0063 | 0.0066 | 0.0875 | Bal. |

### 3.2 Test specimens and fabrication orientations

The details of the lattice designs including the thicknesses, dimensions and relative densities are listed in Table 2, in which 3 groups of the lattice directions of [100], [110] and [111] were fabricated, while a unit cell of $D = 5\ \text{mm}$ is adopted for fabrication based on our preliminary experimental tests to guarantee the fabrication quality. In each group, 6 specimens were built, including 3 variable thickness lattices and 3 uniform thickness ones with the same designed relative density ( $\bar{\rho} = 10.55\%$ ) for comparison. Thus, a total number of 18 specimens were prepared for experimental tests. The weight of the lattices was measured by the dry weight method, and the relative density was calculated according to the volume occupied by the lattices and the density of SS316L, which is also included in Table 2. The low relative standard deviations of the measured relative densities for all groups reveal good fabrication accuracy. The fabrication orientations and specimens built by micro-LPBF along [100], [110] and [111] directions are shown in Fig. 3, in which the lattice directions are labelled by arrows. The material and lattice directions are further illustrated in Table 3, in which the material direction represents the direction in the platform coordinate system shown in Fig. 3. According to the rules of Miller indices, a negative number in the lattice direction is represented by a bar over it. It is noted that the material directions of [100] and [110] specimens were aligned with the lattice directions, while those of [111] specimens were no longer aligned with the lattice directions since the specimens were rotated to save the fabrication space.

To verify the thickness distribution of the as-printed lattices, the uniform and variable thickness [100] specimens were examined by the scanning electron microscope (SEM) using the *JSM-7800F* field emission SEM system. The surface roughness was measured by *RH-2000 High-Resolution 3D Optical Microscope*, while the dimension and thickness variations were evaluated using *Leica DMI6000B Optical Microscope*.

Table 2 Design and measured parameters of uniform and variable thickness N14 shell lattices

| | Uniform thickness | | | Variable thickness | | |
|---|---|---|---|---|---|---|
| Designed thickness $\delta$ (μm) | 122 | | | 100 ~ 250 | | |
| Designed relative density $\bar{\rho}$ | 10.55% | | | 10.55% | | |
| Lattice direction | [100] | [110] | [111] | [100] | [110] | [111] |
| Measured relative density $\bar{\rho}$ | 12.84 ± 0.31 % | 13.03 ± 0.40 % | 13.74 ± 0.63 % | 13.88 ± 0.43 % | 14.11 ± 0.50 % | 12.88 ± 0.34 % |
| Unit cell size $D_1 \times D_2 \times D_3$ (mm) | 5×5×5 | 7.07×7.07 ×5 | 8.66×7.07× 12.25 | 5×5×5 | 7.07×7.07 ×5 | 8.66×7.07× 12.25 |
| No. of unit cells | 5×5×5 | 4×4×5 | 3×4×2.5 | 5×5×5 | 4×4×5 | 3×4×2.5 |
| Total size $L_1 \times L_2 \times L_3$ (mm) | 25×25× 25 | 28.28× 28.28×25 | 25.98×28.28 ×30.62 | 25×25× 25 | 28.28× 28.28×25 | 25.98×28.28 ×30.62 |

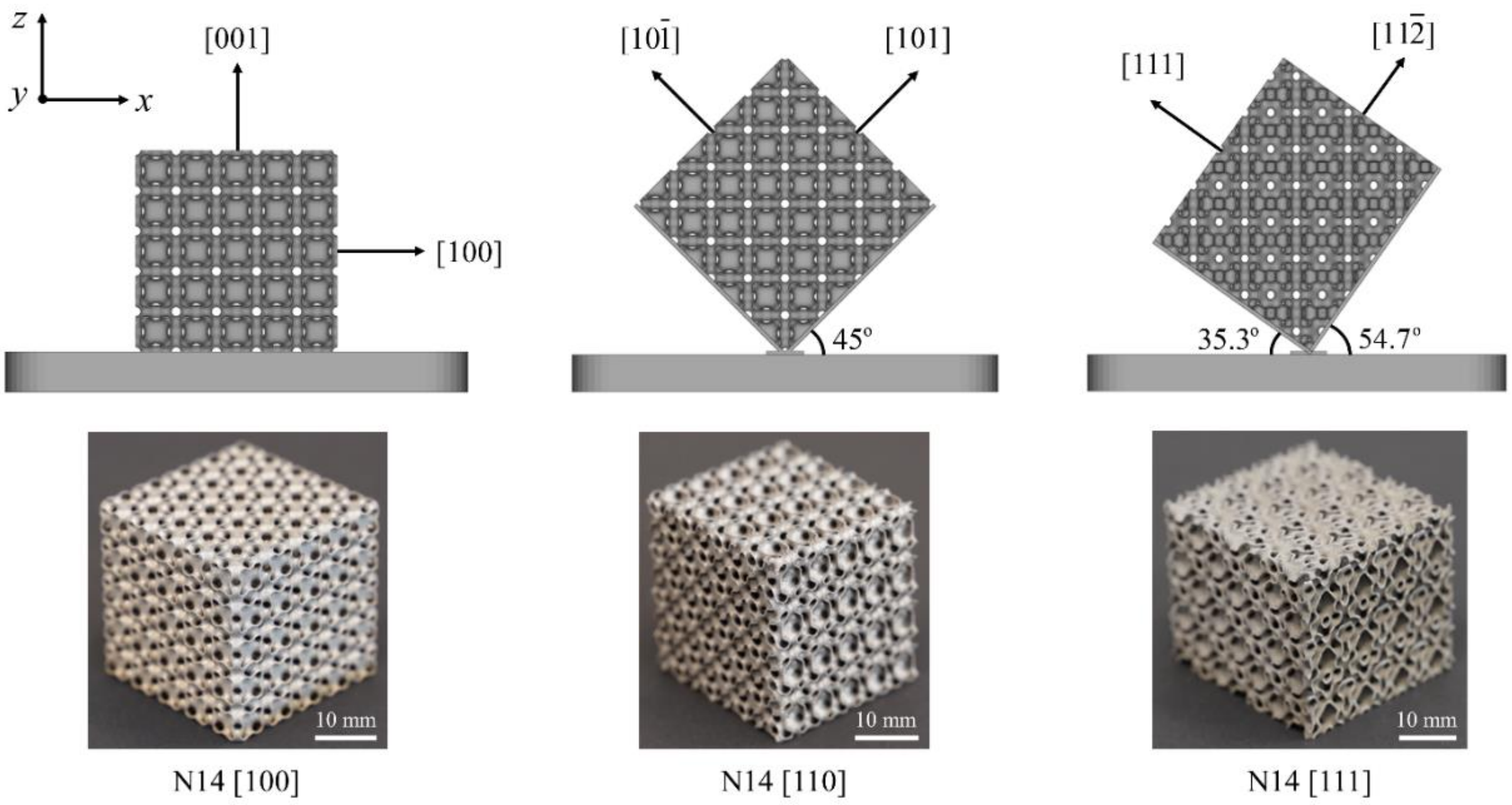

Fig. 3 Fabrication orientations and specimens built by micro-LPBF along [100], [110] and [111] directions. (The lattice directions are labelled by arrows in the figure.)

Table 3 Material and lattice directions of [100], [110] and [111] specimens

| [100] specimen | | [110] specimen | | [111] specimen | |
|---|---|---|---|---|---|
| Material direction | Lattice direction | Material direction | Lattice direction | Material direction | Lattice direction |
| (1, 0, 0) | [100] | (1, 0, 1) | [101] | (0, 1, 0) | $[\bar{1}10]$ |
| (0, 1, 0) | [010] | (1, 0, -1) | $[10\bar{1}]$ | (-0.816, 0, 0.577) | [111] |
| (0, 0, 1) | [001] | (0, 1, 0) | [010] | (0.577, 0, 0.816) | $[11\bar{2}]$ |

### 3.3 Compression tests of elasticity

To test the Young's moduli of the lattices, the *Instron E10000* universal testing system with a 10 kN load cell was used. For brevity of expression, the uniform/variable thickness [100] specimens are denoted as U100/V100, and so forth. In each compression test, the specimen was loaded up to 6 kN with a speed of 0.75 mm/min (corresponding to a strain rate of $4.07\times10^{-4} \sim 5.07\times10^{-4}\,/\,\mathrm{s}$) and unloaded with the same speed. All specimens were within the linear elastic region during the load-unload tests. For each specimen, the compression tests were conducted along all the 3 orthogonal directions. Typically, the compression response of the first loading test deviates from the subsequent ones to some extent due to non-flat surfaces and local yielding, while those of the second and later tests show a good repeatability and the slopes are consistent with the unload slope of the first test, as shown in Fig. 4 for the test of U100 specimen along [100] direction. Thus, the compression test along each direction was conducted for 4 times, with the specimens rotated along the loading axis for $90^{\circ}$ each time. The compressive strain of the lattices was calculated by $\varepsilon = d/h$, where $d$ denotes the crosshead displacement, while $h$ denotes the height of the specimen. The compressive load was measured by the load cell, and the compressive stress was calculated by $\sigma = P/A_0$, where $P$ denotes the compressive load, while $A_0$ denotes the end face area of the specimen. Since the unload curves were more stable than load curves, the Young's modulus of the lattices along each direction was calculated as the average of the 4 unload slopes of the stress-strain curves, denoted as $k_{m+s}$. It should be noted that the deformations of machine and fixtures were included in the crosshead experiments, so the moduli of the lattices were underestimated. To eliminate this underestimation, the machine stiffness was measured by compression tests of the machine itself, following the same load-unload procedure. The compression tests were conducted for 3 times, and the machine stiffness was evaluated as the average of the 3 unload slopes of the load-displacement curves, denoted as $k_m$. Thus, the Young's modulus of the lattices was corrected by:

$$E = \frac{hk_m k_{m+s}}{hk_m - A_0 k_{m+s}}, \tag{19}$$

Eq. (19) was then adopted for evaluation of the Young's moduli of the specimens.

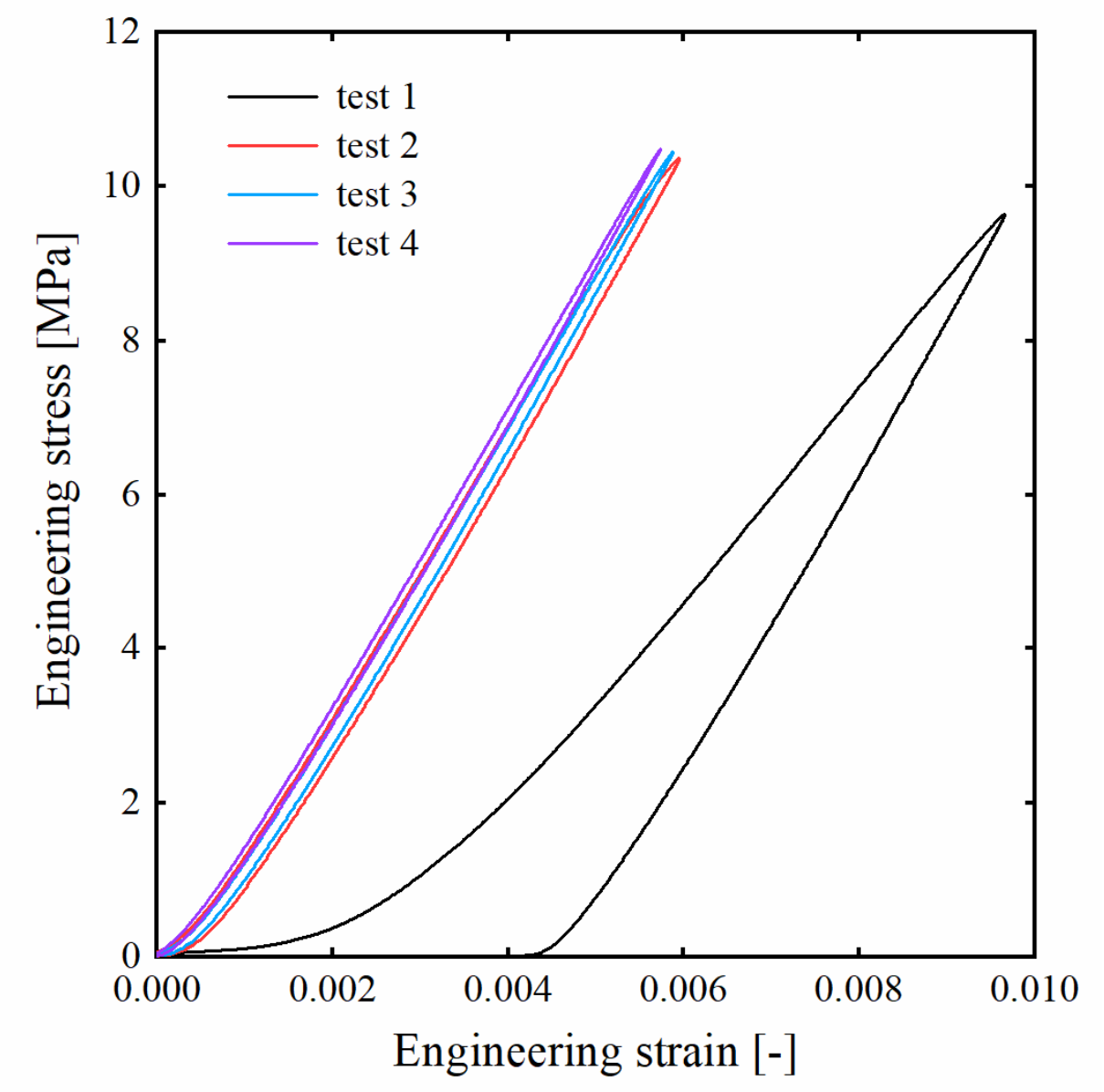

Fig. 4 Stress-strain curves of small strain compression tests along [100] direction for U100 specimen.

### 3.4 Large strain compression tests

To test the failure behavior and energy absorption performance of the lattices, the *MTS ALLIANCE RT/50* testing machine with a 50 kN load cell was used. The tests were conducted along the lattice direction [100] for U100/V100 specimens, [110] for U110/V110 specimens, and [111] for U111/V111 specimens, respectively. In each test, the specimen was loaded up to 45 kN with a compressive speed of 1.5 mm/min (corresponding to a strain rate of $8.89\times10^{-4} \sim 1.01\times10^{-3}\,/\,\mathrm{s}$). The compressive strain and stress were measured from the crosshead displacement and load cell, respectively. The deformations of the lattices were recorded by an *EOS 6D* full frame digital single-lens reflex camera with a frequency of 0.1 Hz. The normalized plateau stress (NPS) $\sigma_p$ was explored as an indicator of the plasticity, which was defined as the average stress within the strain range of [0.2, 0.4] normalized by the relative density, based on ISO 11134:2011 [51]:

$$\sigma_p = \frac{5}{\rho}\int_{0.2}^{0.4} \sigma(e)de, \tag{20}$$

The NPS for the uniform/variable thickness specimens along the three lattice directions [100], [110] and [111] are calculated from experimental data of large strain compression tests by Eq. (20), based on which the maximum to minimum plateau stress ratios are obtained for the uniform/variable thickness specimens, respectively. Besides, the energy absorption performance of the lattices was also studied. Based on the stress-strain curve of the uniaxial compression test, the energy absorption efficiency was calculated as [52]:

$$\eta(\varepsilon) = \frac{1}{\sigma(\varepsilon)}\int_{0}^{\varepsilon} \sigma(e)de, \tag{21}$$

The strain with the maximum efficiency was then extracted, denoted as $\varepsilon_d$, based on which the specific energy absorption (SEA) was evaluated as:

$$\psi = \frac{1}{\rho}\int_{0}^{\varepsilon_d} \sigma(e)de, \tag{22}$$

to measure the energy absorption capacity of the lattices.

# 4. Results and discussions

## 4.1 Elastically-isotropic designs

The six types of TPMS in Fig. 1 are adopted for achieving elastically-isotropic shell lattices based on the optimization problem (17). The unit cell size is set as $D = 2$ mm, while the initial design is taken as the uniform thickness lattice with the thickness of $\delta = 0.04$ mm, corresponding to an aspect ratio $\delta / D$ (i.e., shell thickness to unit cell size ratio) of 1/50. The upper bound is set as $(1/10)D$ for the thin shell assumption to hold, while the lower bound should be properly selected for the optimization algorithm to yield elastically-isotropic designs satisfying the manufacturing constraints. By adopting the minimum thickness $\delta_{\min}$ of $(1/50)D$ and maximum thickness $\delta_{\max}$ of $(1/10)D$, the elastically-isotropic N14 lattice with the relative density of 10.55% is achieved, which is used for experimental tests. For the targeted relative density of $\rho^* = 10\%$, the appropriate values of the minimum thicknesses are found to be $(1/60)D$, $(1/60)D$, $(1/100)D$, $(1/150)D$, $(1/250)D$ and $(1/1200)D$ for N14, OCTO, IWP, FRD, N and P lattices, respectively. The optimized elastically-isotropic variable thickness TPMS lattices with the relative density of 10% are shown in Fig. 5, with the thickness distributions and views of the 1/48 fundamental unit, 1/8 cell and unit cell, respectively, while the geometry and mechanical data of the uniform and variable thickness lattices are listed in Table 4, in which the average mesh size to unit cell size ratio $\Delta / D$ is included and shows that the fine mesh was used in the FEA and optimization process. The stiffness properties of the uniform thickness lattices in Table 4 refer to those along [100] direction, while the lattices possess different stiffness properties along different lattice directions due to the elastic anisotropy.

Given the relative density of $\overline{\rho} = 10\%$, the HS upper bounds for the normalized Young's modulus $E / (\overline{\rho} E_s)$, shear modulus $G / (\overline{\rho} E_s)$ and bulk modulus $K / (\overline{\rho} E_s)$ in Eq. (10a) ~ (10b) are derived as 0.525, 0.212 and 0.338, respectively. All the uniform thickness TPMS shell lattices are found to approach nearly 94% of the HS upper bound of bulk modulus in such relative density. The uniform thickness N14/OCTO shell lattices possess lower elastic anisotropy than the remaining lattices, making it easier to achieve elastic isotropy by varying the shell thickness, with the Young's, shear and bulk moduli reaching 56.57%/55.81%, 52.83%/52.36% and 90.83%/84.62% of the HS upper bounds, respectively. In contrast, the higher elastic anisotropy of the uniform thickness IWP/FRD/N lattices makes it harder to achieve elastic isotropy by varying the shell thickness, while smaller lower bounds $\delta_{\min}$ are required for the optimization algorithm to yield elastically-isotropic designs, with the Young's, shear and bulk moduli reaching 48.57%/43.24%/43.24%, 46.23%/40.09%/40.57% and 64.50%/67.75%/63.02% of the HS upper bounds, respectively. Thus, the bulk moduli of the elastically-isotropic IWP/FRD/N lattices are lower than those of N14/OCTO, as a result of the larger deviation from uniform thickness. The uniform thickness P lattice possesses the highest elastic anisotropy among the six TPMS shell lattices, and a significantly larger design space is needed for elastically-isotropic designs. The optimized elastically-isotropic P lattice therefore possesses the widest range of thickness distribution, with a maximum to minimum thickness ratio of nearly 120, while the Young's, shear and bulk moduli are the lowest among the six types of elastically-isotropic TPMS lattices, reaching 26.48%, 24.53% and 51.48% of the HS upper bounds, respectively. In Callens et al.'s work [36] that achieved elastic isotropy by two constituent materials, a large Young's moduli ratio of 100 was used to achieve elastically-isotropic P lattices.

To illustrate the influence of thickness variations on the stiffness properties of the lattices, the normalized Young's moduli along all directions of the six types of uniform and variable thickness lattices in Table 4 are compared in Fig. 6, in which the sphere represents those of the elastically-isotropic variable thickness lattices, while the other transparent surface represents those of the elastically-anisotropic uniform ones. The results show that the spheres of elastically-isotropic N14, OCTO and N lattices intersect with the surfaces of the uniform ones, revealing that the Young's moduli of the elastically-isotropic lattices are in between the maximum and minimum values of those of the uniform ones with equal relative densities. On the other hand, the spheres of the elastically-isotropic IWP, FRD and P lattices are fully contained within the surfaces of the uniform thickness lattices, meaning the stiffness of the optimized lattices is decreased to reach elastic isotropy. Therefore, the increase or decrease of the Young's moduli of the optimized lattices compared with the uniform ones depends on the type of TPMS. For N14, OCTO and N lattices, the optimization procedure enhances the weaker directions and decreases the stiffer directions, leading to a balanced stiffness. For IWP, FRD and P lattices, the elastic properties along all directions are decreased, thus leading to lower stiffness.

Furthermore, the optimization workflow is implemented to obtain the six types of elastically-isotropic TPMS lattices with relative densities of 5%, 15% and 20%, respectively, in which the upper and lower bounds of the shell thicknesses are adjusted accordingly. The normalized Young's and bulk moduli of the optimized lattices are plotted versus the relative density in Fig. 7, which reveals that the stiffness properties of all elastically-isotropic TPMS lattices (except for P) increase slightly as the relative densities increase, while their bulk moduli are closer to the HS upper bounds than the Young's moduli. Among the six types of elastically-isotropic TPMS lattices, N14 and OCTO exhibit the most superior stiffness properties, while their Young's and bulk moduli are always kept over 50% and 80% of the HS upper bounds in low relative densities, respectively.

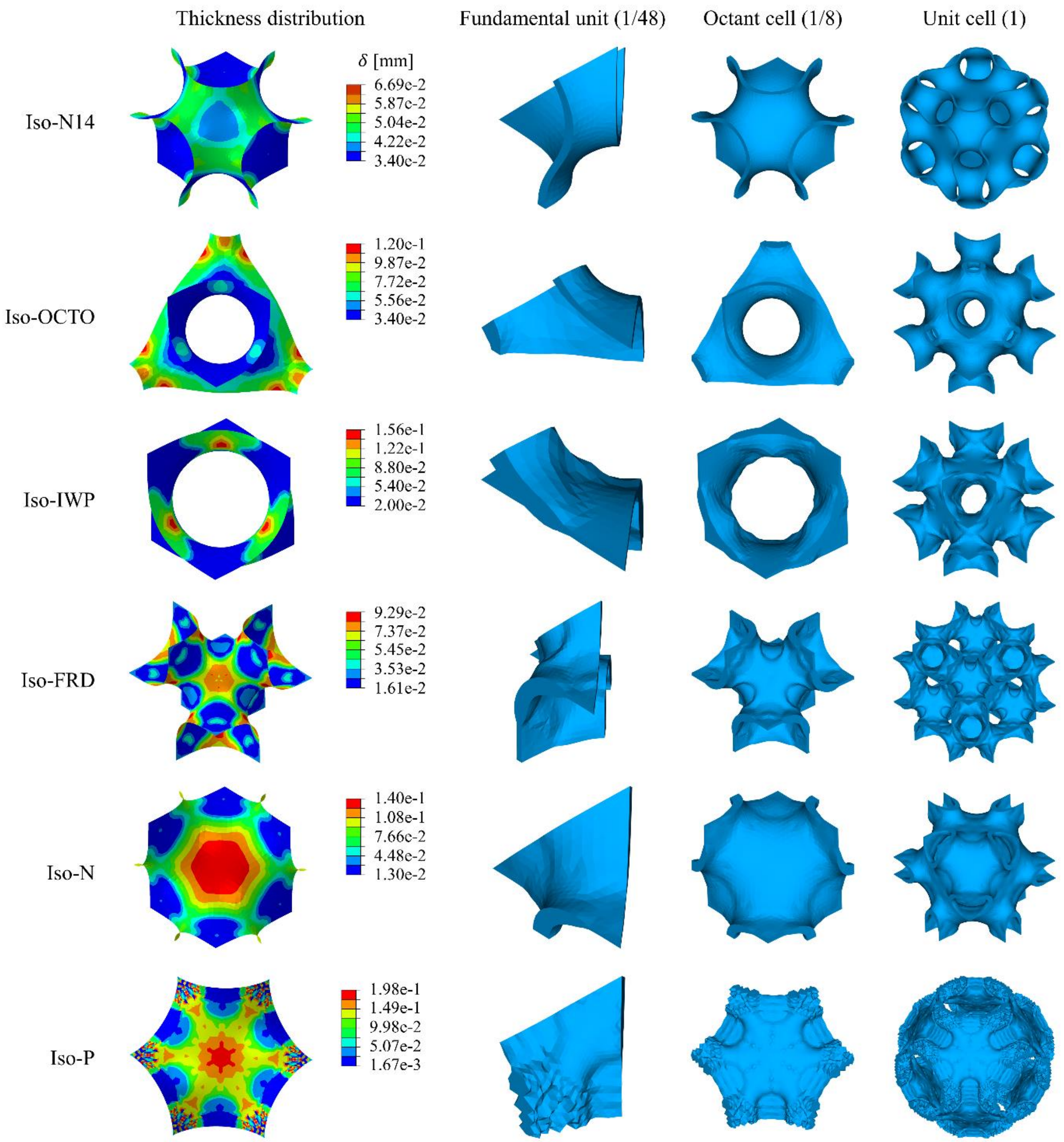

Fig. 5 Elastically-isotropic variable thickness TPMS shell lattices with the relative density of $\bar{\rho} = 10\%$: thickness distributions and views of the 1/48 fundamental unit, 1/8 cell and unit cell. (*Iso-N14* represents elastically-isotropic N14 shell lattices, and so forth.)

Table 4 Geometry and mechanical data of uniform and variable thickness TPMS shell lattices

| Property | | $\Delta / D$ | $\bar{\rho}$ | $\xi$ | $\nu$ | $E/(\bar{\rho}E_s)$ | $G/(\bar{\rho}E_s)$ | $K/(\bar{\rho}E_s)$ |
|---|---|---|---|---|---|---|---|---|
| Uniform thickness | N14 | 0.017 | 9.91% | 0.839 | 0.321 | 0.341 | 0.108 | 0.317 |
| | OCTO | 0.027 | 10.00% | 0.665 | 0.301 | 0.379 | 0.097 | 0.317 |
| | IWP | 0.021 | 10.00% | 0.609 | 0.271 | 0.436 | 0.104 | 0.317 |
| | FRD | 0.011 | 9.99% | 0.535 | 0.279 | 0.421 | 0.088 | 0.317 |
| | N | 0.021 | 9.97% | 0.398 | 0.283 | 0.414 | 0.064 | 0.317 |
| | P | 0.010 | 10.09% | 2.884 | 0.412 | 0.168 | 0.172 | 0.317 |
| Variable thickness | N14 | 0.017 | 9.91% | 1.009 | 0.339 | 0.297 | 0.112 | 0.307 |
| | OCTO | 0.027 | 10.00% | 1.009 | 0.329 | 0.293 | 0.111 | 0.286 |
| | IWP | 0.021 | 10.00% | 0.999 | 0.305 | 0.255 | 0.098 | 0.218 |
| | FRD | 0.011 | 9.99% | 1.003 | 0.335 | 0.227 | 0.085 | 0.229 |
| | N | 0.021 | 9.97% | 0.997 | 0.322 | 0.227 | 0.086 | 0.213 |
| | P | 0.010 | 10.09% | 1.009 | 0.367 | 0.139 | 0.052 | 0.174 |

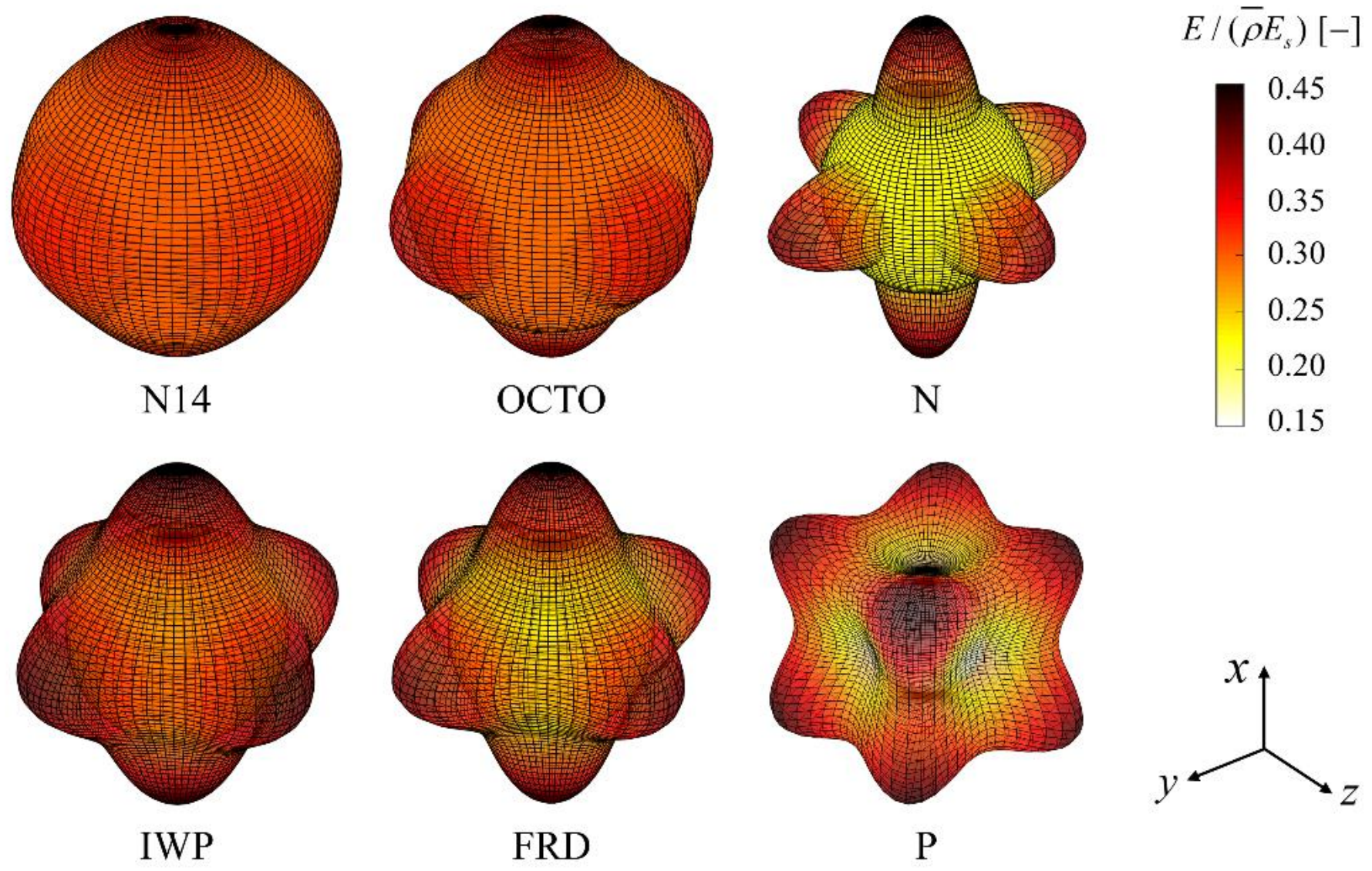

Fig. 6 The normalized Young's moduli for the six types of uniform and variable thickness TPMS shell lattices with the relative density of 10%.

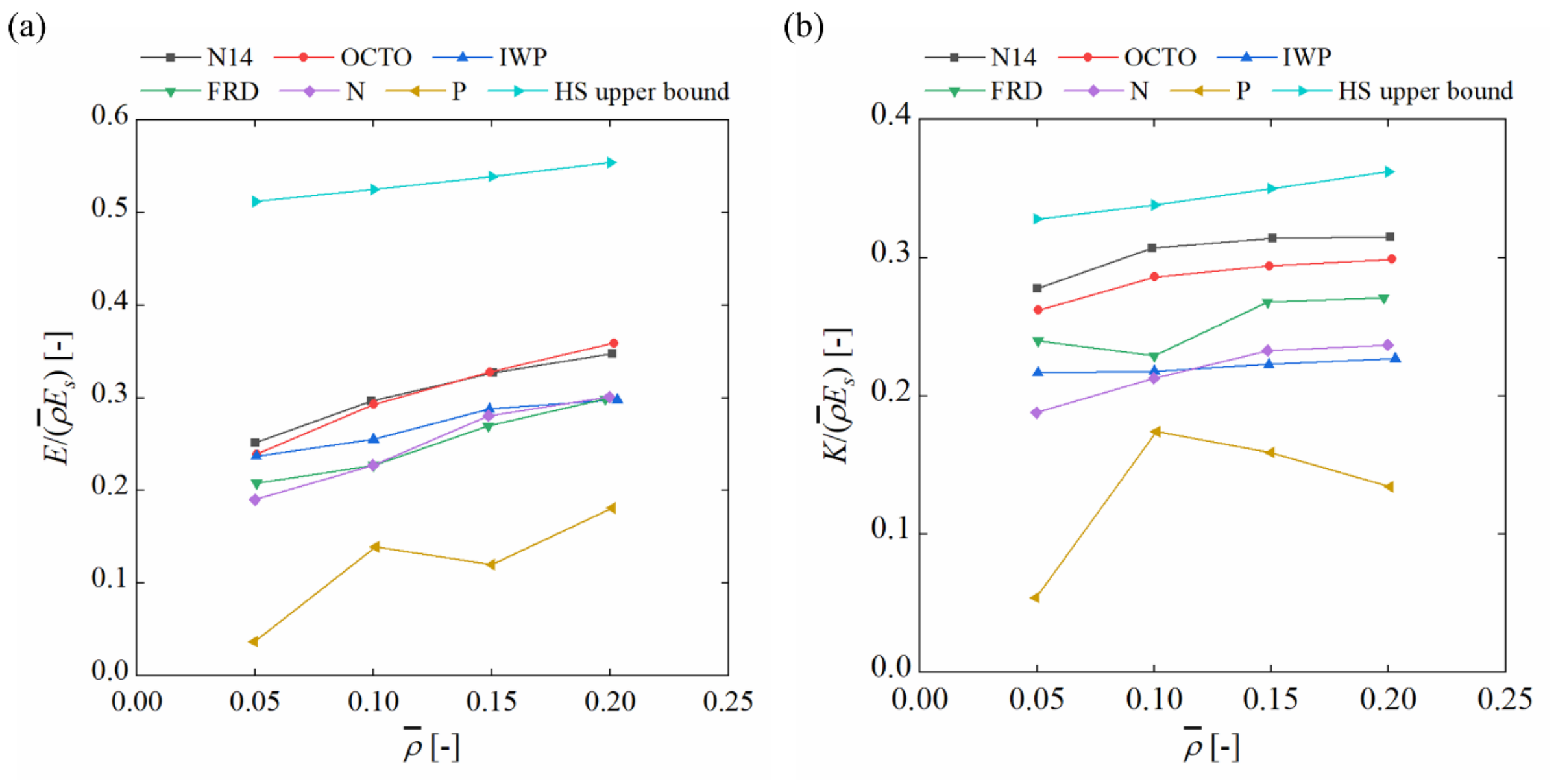

Fig. 7 (a) Normalized Young's modulus, (b) normalized bulk modulus versus the relative density, for elastically-isotropic TPMS shell lattices.

### 4.2 Stiffness and initial yield strength comparisons

The stiffness properties of the six types of elastically-isotropic TPMS shell lattices are compared with the elastically-isotropic plate lattices and truss lattices with the same relative density of $\bar{\rho} = 10\%$, as summarized in Fig. 8, in which Opt Iso-PLS and Opt Iso-TLS represent the optimal elastically-isotropic plate lattices and truss lattices [12-17], respectively, while Iso-hTLS represents the elastically-isotropic hollow truss lattices and includes three elementary classes: SC, body-centered cubic (BCC) and face-centered cubic (FCC) [53]. The stiffness properties of the optimal elastically-isotropic plate lattices are validated to approach the HS upper bounds in low relative density limits [15-17], while those of the optimal elastically-isotropic truss lattices can only reach nearly 1/3 of the HS upper bounds in low relative density limits [12-14]. Given the relative density of 10%, the stiffness properties of the optimal elastically-isotropic plate lattices are numerically confirmed to approach over 98% of the HS upper bounds [16], while those of the optimal elastically-isotropic truss lattices are found to reach nearly 40% of the HS upper bounds [14, 53], as shown in the red and black crossed squares in Fig. 8, respectively. The grey region in Fig. 8 represents the stiffness properties achievable by elastically-isotropic truss lattices, in which the three classes of elastically-isotropic hollow truss lattices possess lower stiffness properties than the optimal elastically-isotropic truss lattices.

In contrast, the bulk moduli of the elastically-isotropic N14/OCTO shell lattices are found to approach nearly 90%/85% of the HS upper bounds, those of the elastically-isotropic IWP/FRD/N lattices can reach over 60% of the HS upper bounds, while those of the elastically-isotropic P lattices can reach nearly 50% of the HS upper bounds. On the other hand, the Young's and shear moduli of the elastically-isotropic N14/OCTO/IWP/FRD/N lattices reach nearly 40% ~ 60% of the HS upper bounds, while those of the elastically-isotropic P lattices can only reach 20% ~ 30% of the HS upper bounds. Therefore, the majority of elastically-isotropic shell lattices exhibit superior stiffness properties and significantly outperform elastically-isotropic truss lattices of equal relative densities, while their stiffness properties are lower than the elastically-isotropic plate lattices, as a trade-off induced by the open-cell topology. In general, the elastically-isotropic plate lattices maximize the stiffness properties of porous cellular solids in low relative densities, while their closed-cell topology complicates the manufacturing and applications in specified industries. The elastically-isotropic shell lattices developed in this work exhibit superior stiffness properties and maintain the open-cell topology. The opportunities for elastically-isotropic shell lattices are shown in the red region in Fig. 8, while the detailed upper limits of their stiffness properties remain to be explored. Thus, it is an ongoing research topic to develop elastically-isotropic open-cell shell lattices with the stiffness properties closer to the HS upper bounds.

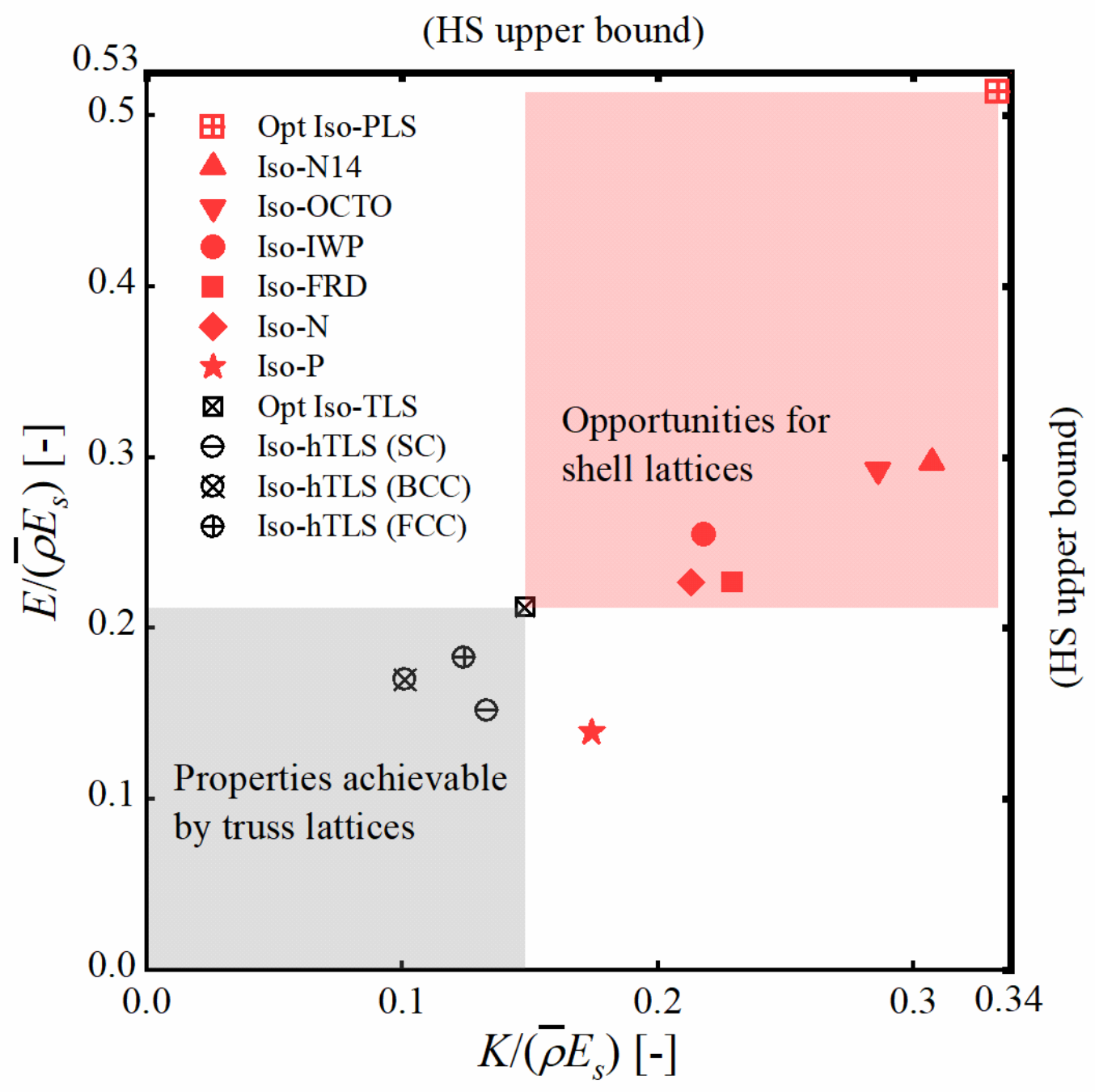

Fig. 8 Comparisons of the stiffness properties of elastically-isotropic plate lattices, shell lattices and truss lattices, presented in the form of normalized Young's modulus versus normalized bulk modulus.

To illustrate the influence of thickness variations on the initial yield stresses, the von Mises stress contours of the uniform and variable thickness N14 shell lattices under the uniaxial stress loading, hydrostatic loading and pure shear loading are shown in Fig. 9. The relative density of both types of lattices is $\bar{\rho} = 10.55\%$, while the corresponding normalized initial yield stresses are listed in Table 5. The thickness variations are shown to decrease the normalized initial yield stresses by 12.82% and 18.35% for the uniaxial stress loading and hydrostatic loading, respectively, while increase that of the pure shear loading by 21.54%.

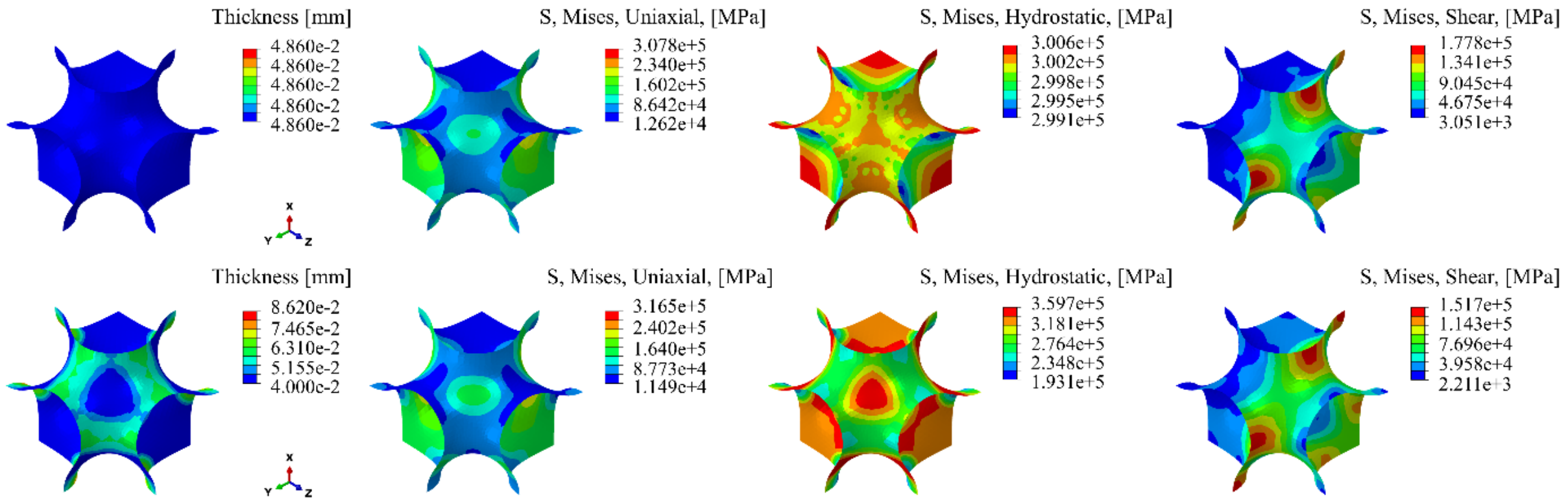

Fig. 9 Thickness distributions and von Mises stress contours of the uniform and variable thickness N14 shell lattices. The first to fourth columns represent the thickness distributions, von Mises stress contours under the uniaxial stress, hydrostatic and pure shear loading, respectively.

Table 5 Comparisons of the normalized initial yield stresses of the uniform and variable thickness N14 shell lattices under the three load cases

| | $\bar{\rho}$ | $(\frac{\Sigma_s^y}{\rho\sigma_s^y})^{(US)}$ | $(\frac{\Sigma_s^y}{\rho\sigma_s^y})^{(T)}$ | $(\frac{\Sigma_s^y}{\rho\sigma_s^y})^{(S)}$ |
|---|---|---|---|---|
| Uniform thickness | 10.55% | 0.234 | 0.665 | 0.130 |
| Variable thickness | 10.55% | 0.204 | 0.543 | 0.158 |

### 4.3 Morphologies of fabricated lattice specimens

The SEM images of the top surfaces are shown in Fig. 10, which demonstrates the uniform thickness distribution in Fig. 10 (a) and variable thickness distribution in Fig. 10 (b). Further observations of the molten pool of the micro-LPBF fabricated lattices from the top surface reveal a pool size of approximately 50 μm. Therefore, the micro/high-precision LPBF configuration results in higher dimensional accuracy, smaller surface roughness and part distortion as compared with the conventional LPBF machine [47].

Furthermore, the formation of bonded powders at the hanging features facing the building platform inevitably have effects on the geometry deviations and surface roughness [54, 55], which influences the mechanical properties of the fabricated specimens. The surface roughness $Ra$ is shown in Fig. 11, implying that the top surface possesses the lowest roughness with an $Ra$ of 3.18 μm, while the bottom surface has the largest roughness with an $Ra$ of 12.93 μm due to the sintered powders. The dimension and thickness variations of the specimens are shown in Fig. 12, revealing that the dimension variations are all below 5%, while the thickness variations are below 10% and 15% for the top and lateral surfaces, respectively. The fabricated specimens are therefore found to be quite consistent with the designed models, with minor geometry deviations.

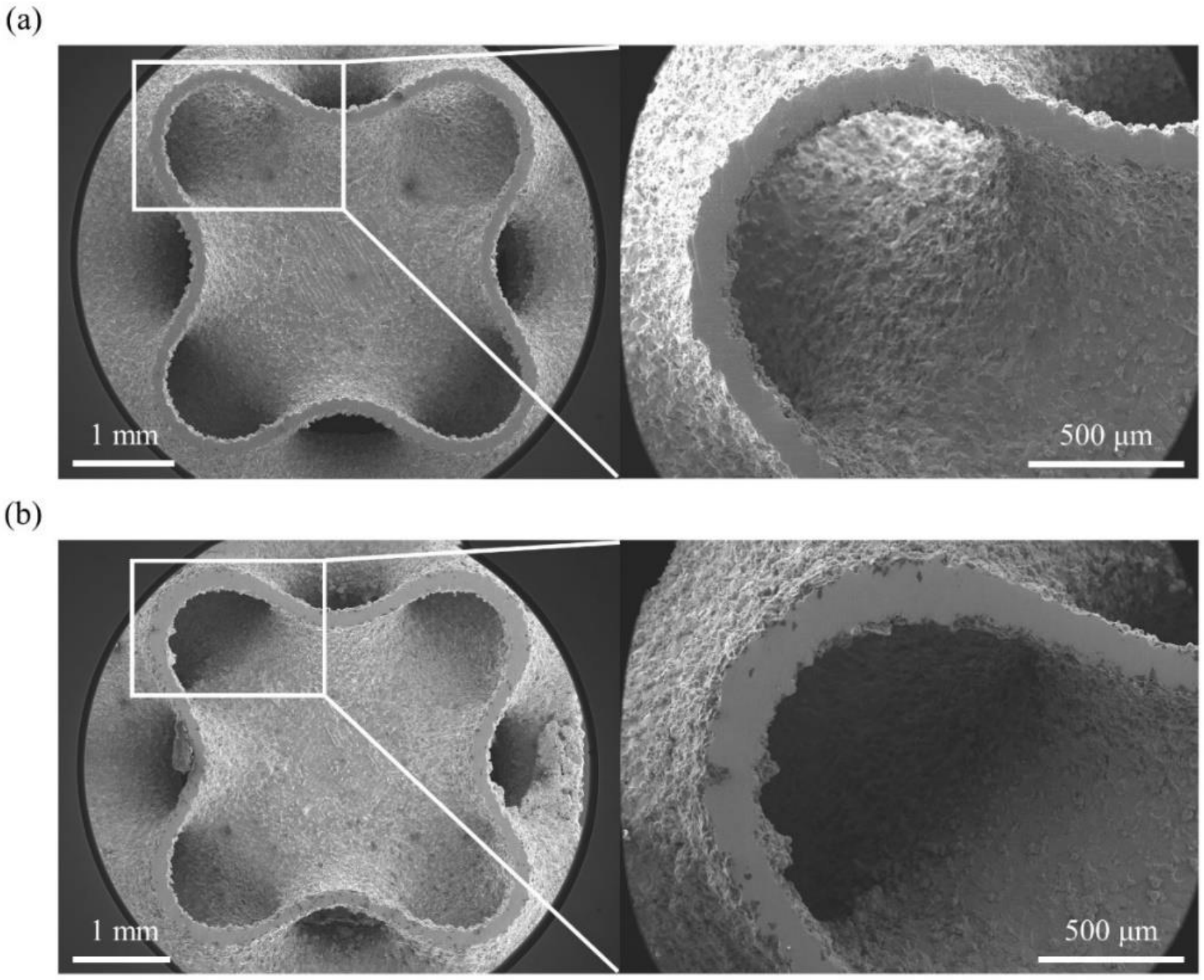

Fig. 10 SEM images of [100] specimens: (a) uniform thickness lattices, (b) variable thickness lattices.

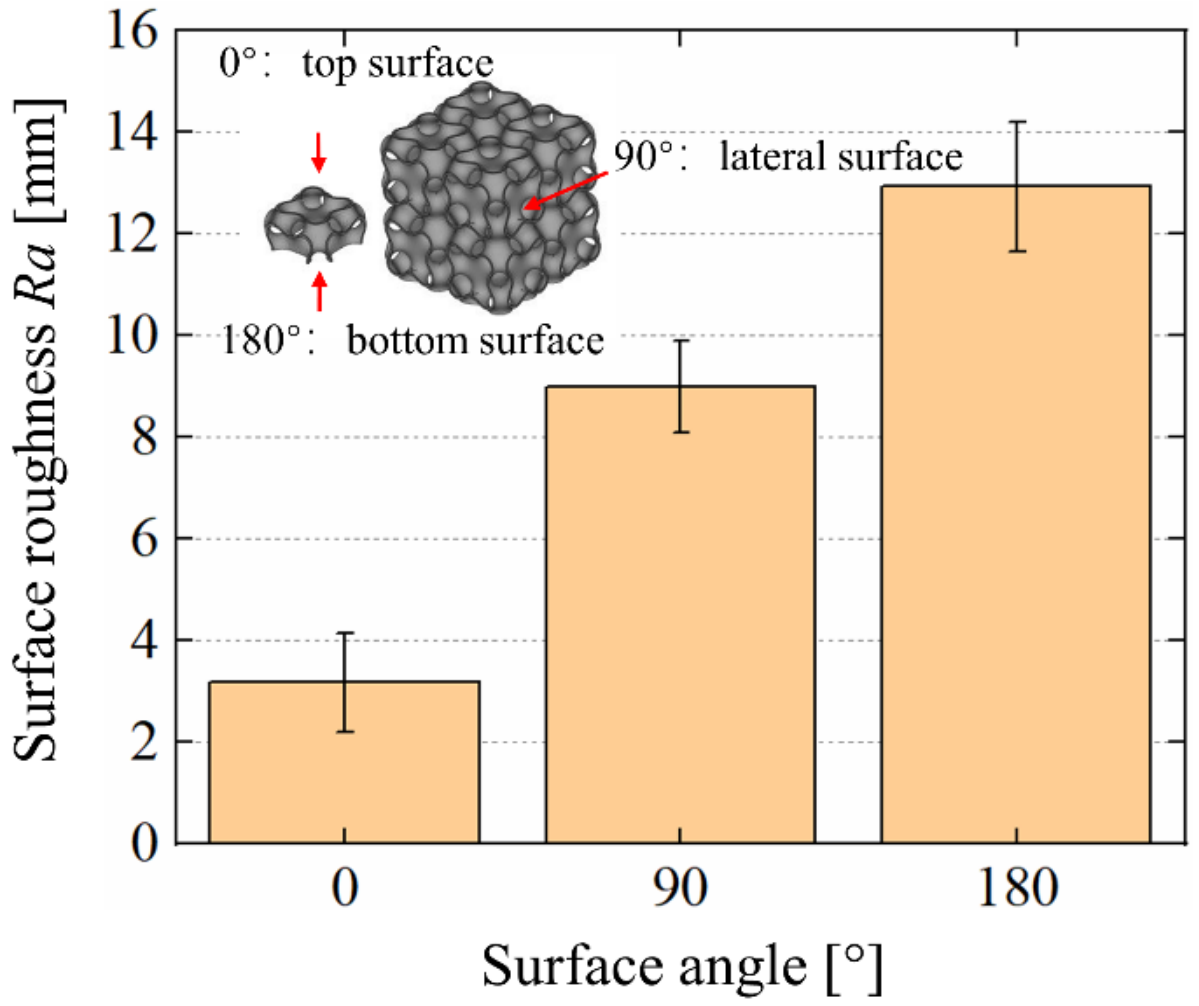

Fig. 11 The surface roughness measured for various angles.

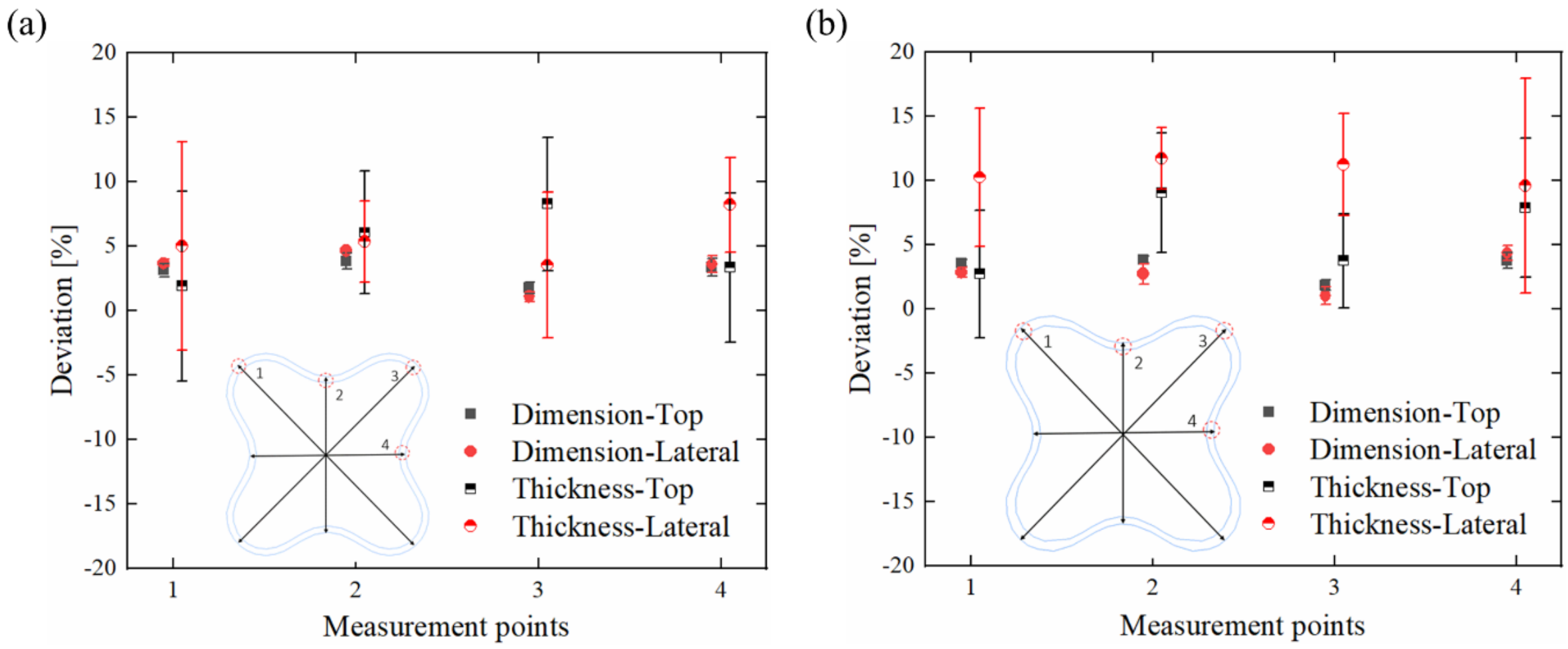

Fig. 12 Dimension and thickness deviations of the fabricated N14 lattice specimens: (a) uniform thickness, (b) variable thickness.

### 4.4 Experimental results of elasticity

The experimental and FEA results of the normalized Young's moduli of the uniform and variable thickness specimens are listed in Table 6. It is found that the multiple load-unload tests lead to high experimental repeatability. Theoretically, the Young's moduli along [100], [010] and [001] directions are the same, while comparisons of the experimental data along (1), (2) and (3) directions reveal that the Young's moduli along [100] and [010] directions are almost the same, which, however, are nearly 13% higher than that along [001] direction. Observations of the Young's moduli along (4), (5) and (7) directions reveal a similar trend. Thus, the elastic anisotropy of the constituent material is observed and yields a detectable deviation of the stiffness properties of the AM fabricated lattices. Since [001] is the building direction, while [100] and [010] lie in the platform ($x$ - $y$) plane, two key observations are made that: (1) the LPBF fabricated material exhibits similar properties within the platform plane, which are different from the building direction ($z$ direction); (2) the larger the angle between the specified direction and the platform plane, the lower the Young's modulus induced by the constituent material's elastic anisotropy.

The normalized Young's moduli along all the 9 test directions for the uniform and variable thickness specimens are further illustrated in Fig. 13. The 9 tests are categorized into 5 groups, denoted as GR1-GR5, each of which has the same material direction and therefore the same material property along the test direction. For the test directions (1), (2), (6) and (7) within GR1, the material directions all lie in the platform plane. Since (1), (2) and (6) also possess the same lattice direction of <100>, the tested normalized Young's moduli are almost the same, with an average of 0.275/0.261 for the uniform/variable thickness specimens. However, (7) has a different lattice direction of <110>, thus with the different normalized Young's moduli of 0.254/0.255. Therefore, the experimental data reveals a [100] to [110] moduli ratio of 1.08/1.02 for the uniform/variable thickness specimens, respectively. By comparison, the FEA gives the [100] to [110] moduli ratio as 1.12/1.00, which shows a reasonable agreement with the experimental results. In addition, the FEA gives the normalized [100] Young's moduli of 0.343/0.306 for the uniform/variable thickness specimens, showing 24.73%/17.24% relative differences with the experimental results; and normalized [110] Young's moduli of 0.306 for both specimens, showing 20.47%/20.00% relative differences. As a whole, the differences between experimental and FEA results are within an acceptable range, while the remaining discrepancy is attributed to: (1) the boundary conditions of multi-cell specimens in experiments are approximations to the periodic boundary conditions in FEA; (2) manufacturing imperfections such as surface roughness; (3) the elastic anisotropy of the constituent material.

Besides, the experimental [100] to [110] moduli ratio reveals a reduced elastic anisotropy of the variable thickness specimens, with the ratio reduced from 1.08 of the uniform thickness specimens to 1.02 of the variable thickness ones. For the test direction (8) with the lattice direction of [111], and (9) with that of [$11\bar{2}$], the material directions do not lie in the platform plane and has the intersection angles of $35.3^{\circ}$ and $54.7^{\circ}$, respectively. Thus, the tested Young's moduli deviate more from the FEA results, as a result of the elastic anisotropy of the AM fabricated material. Overall, the maximum to minimum Young's moduli ratio of all the 9 test directions is reduced from 1.26 of the uniform thickness specimens to 1.20 of the variable thickness specimens. Although the overall anisotropy in stiffness is still observed in the variable thickness specimens due to the constituent material's anisotropy, it is validated that the elastic anisotropy of the lattices is reduced by varying the shell thickness accordingly.

Table 6 Normalized Young's modulus of uniform and variable thickness N14 shell lattices

| Specimen | Test direction | Material direction | Lattice direction | $E/(\bar{\rho}E_s)$ | |
|---|---|---|---|---|---|
| | | | | Experimental | FEA |
| U100 | (1) | (1, 0, 0) | [100] | $0.279 \pm 0.009$ | 0.343 |
| | (2) | (0, 1, 0) | [010] | $0.276 \pm 0.001$ | |
| | (3) | (0, 0, 1) | [001] | $0.245 \pm 0.005$ | |
| U110 | (4) | (1, 0, 1) | [101] | $0.224 \pm 0.006$ | 0.306 |
| | (5) | (1, 0, -1) | $[10\bar{1}]$ | $0.222 \pm 0.006$ | 0.306 |
| | (6) | (0, 1, 0) | [010] | $0.269 \pm 0.025$ | 0.343 |
| U111 | (7) | (0, 1, 0) | $[1\bar{1}0]$ | $0.254 \pm 0.017$ | 0.306 |
| | (8) | (-0.816, 0, 0.577) | [111] | $0.251 \pm 0.013$ | 0.295 |
| | (9) | (0.577, 0, 0.816) | $[11\bar{2}]$ | $0.227 \pm 0.016$ | 0.306 |
| V100 | (1) | (1, 0, 0) | [100] | $0.259 \pm 0.002$ | 0.306 |
| | (2) | (0, 1, 0) | [010] | $0.257 \pm 0.011$ | |
| | (3) | (0, 0, 1) | [001] | $0.228 \pm 0.006$ | |
| V110 | (4) | (1, 0, 1) | [101] | $0.222 \pm 0.009$ | 0.306 |
| | (5) | (1, 0, -1) | $[10\bar{1}]$ | $0.221 \pm 0.003$ | |
| | (6) | (0, 1, 0) | [010] | $0.266 \pm 0.003$ | |
| V111 | (7) | (0, 1, 0) | $[1\bar{1}0]$ | $0.255 \pm 0.011$ | 0.306 |
| | (8) | (-0.816, 0, 0.577) | [111] | $0.238 \pm 0.028$ | |
| | (9) | (0.577, 0, 0.816) | $[11\bar{2}]$ | $0.222 \pm 0.009$ | |

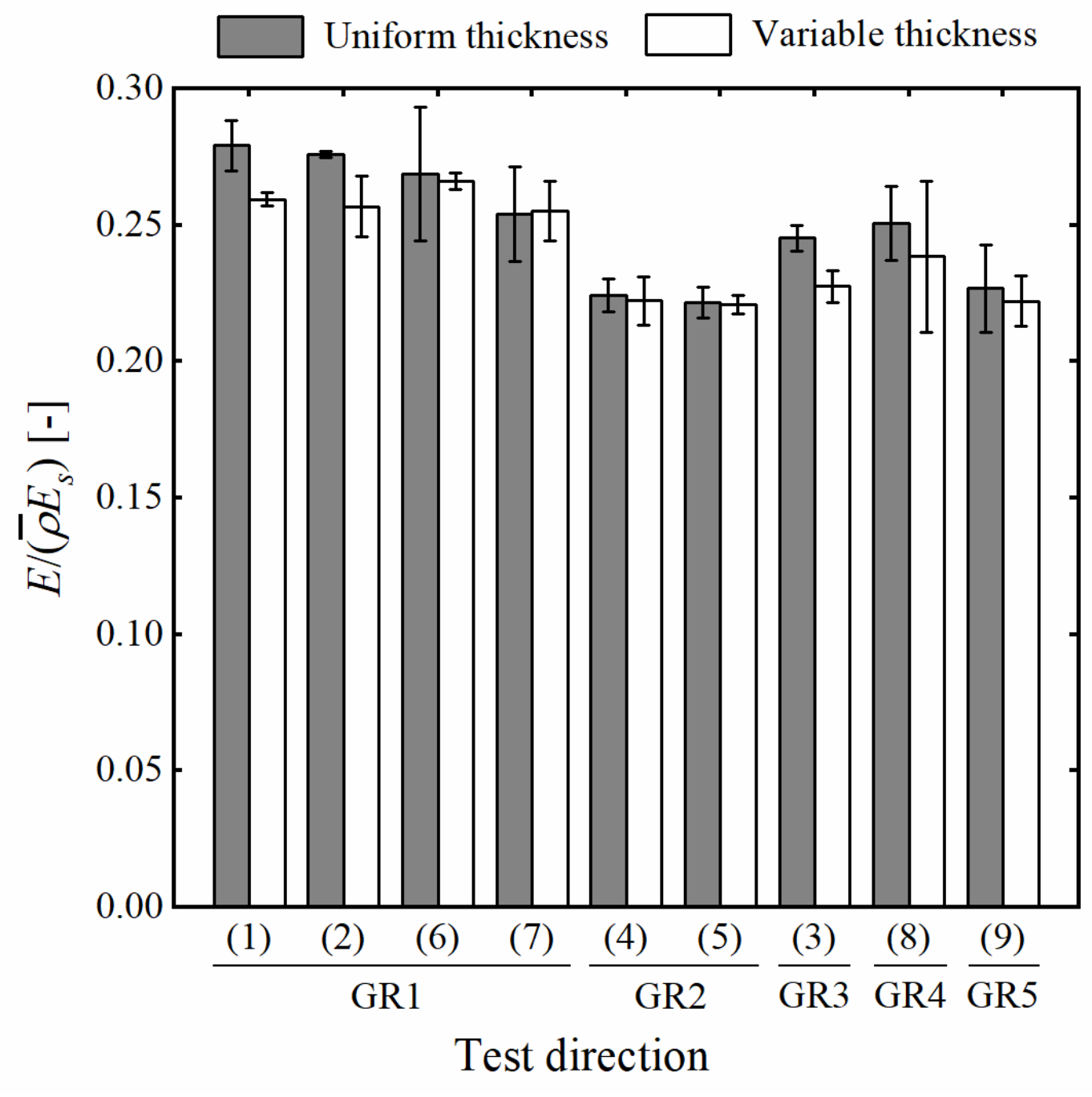

Fig. 13 Normalized Young's moduli along all the 9 test directions for the uniform and variable thickness specimens.

### 4.5 Experimental results of large strain compressions

The stress-strain curves of large strain compression tests along the 3 lattice directions are shown in Fig. 14, with the deformation patterns of the uniform thickness specimens shown in Fig. 15. Experimental results reveal that both uniform and variable thickness N14 lattices show different failure deformation behaviors along different directions. The thickness variations in this work have no significant influences on the failure deformations. It is noted that a layer-by-layer plastic buckling failure mode is observed in the compression tests along [100] direction in Fig. 15 (a), corresponding to 5 sudden drops in the stress-strain curves in Fig. 14 (a). However, no buckling deformation occurs in the compression tests along [110] and [111] directions, while the failures are related to the shear deformation. Thus, no sudden drops appear in the stress-strain curves in Fig. 14 (b) and Fig. 14 (c).

The NPS $\sigma_p$ and SEA $\psi$ are illustrated in Fig. 16, showing that the uniform thickness specimens possess the highest and lowest plateau stresses along [111] and [100] directions, respectively. The variable thickness ones possess the highest and lowest plateau stresses along [110] and [100] directions, respectively. The maximum to minimum plateau stress ratio is reduced from 1.21 of the uniform thickness specimens to 1.12 of the variable thickness ones, which reveals a reduced anisotropy of strength for the variable thickness specimens. Besides, observations of the SEA along the 3 directions reveal that the uniform thickness specimens possess the highest and lowest energy absorption capacities along [111] and [110] directions, respectively. The maximum to minimum SEA ratio is reduced from 1.23 of the uniform thickness specimens to 1.03 of the variable thickness ones, which indicates that the variable thickness specimens can even achieve nearly isotropic energy absorption capacity.

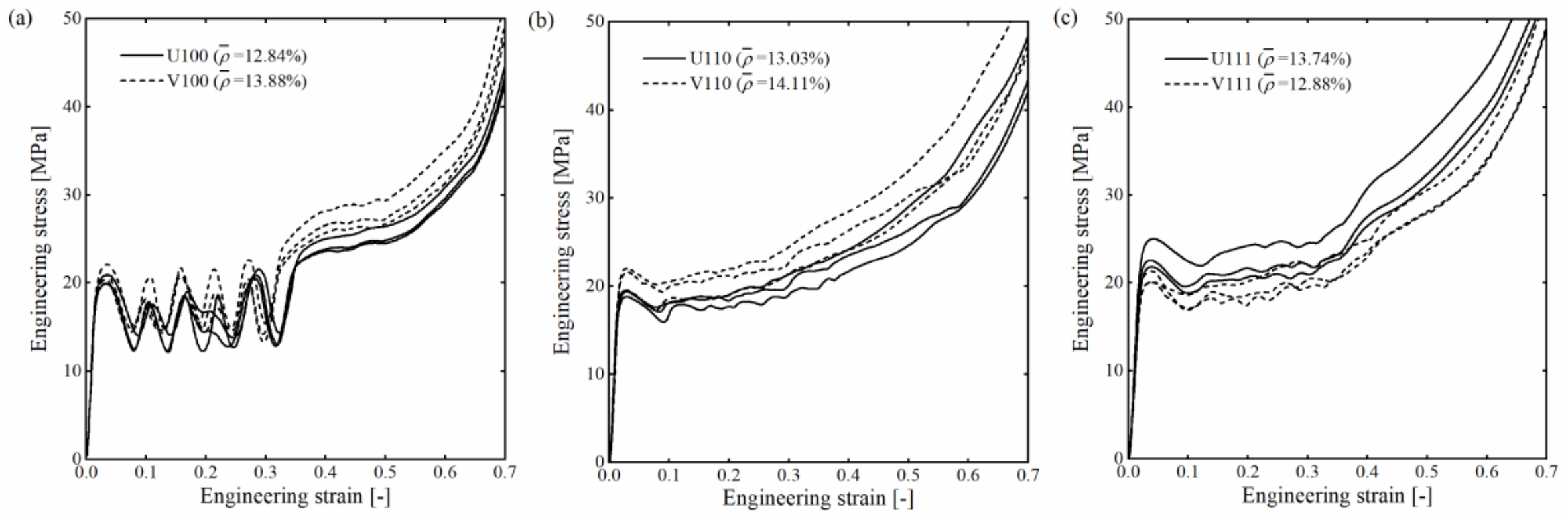

Fig. 14 Stress-strain curves of large strain compression tests along: (a) [100], (b) [110], (c) [111] directions.

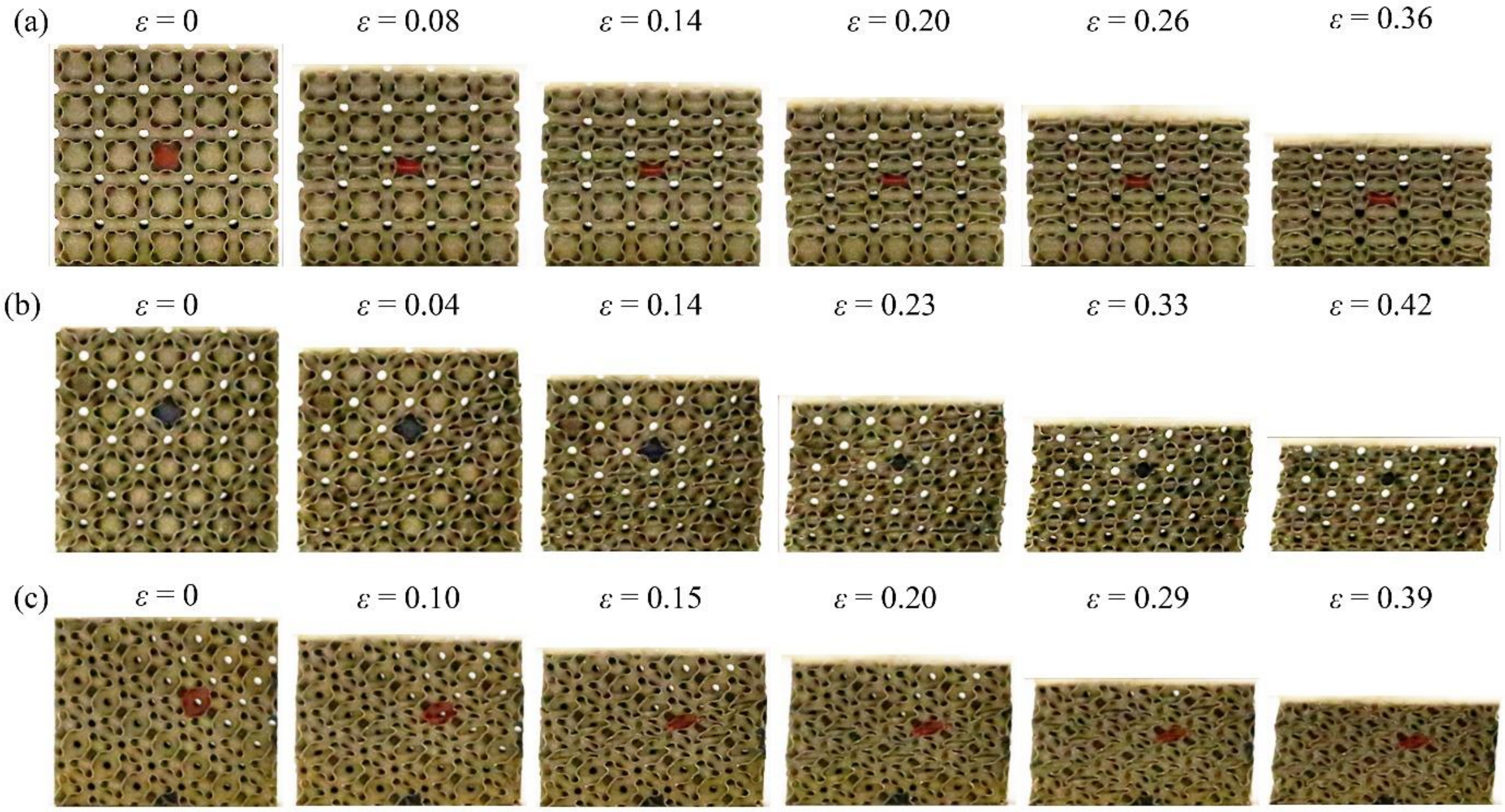

Fig. 15 Experimental large strain compressive deformation patterns for: (a) U100 specimen with a layer-by-layer plastic buckling failure mode, (b) U110 specimen with shear deformation failure, (c) U111 specimen with shear deformation failure. (The corresponding engineering strains are labelled in the figure.)

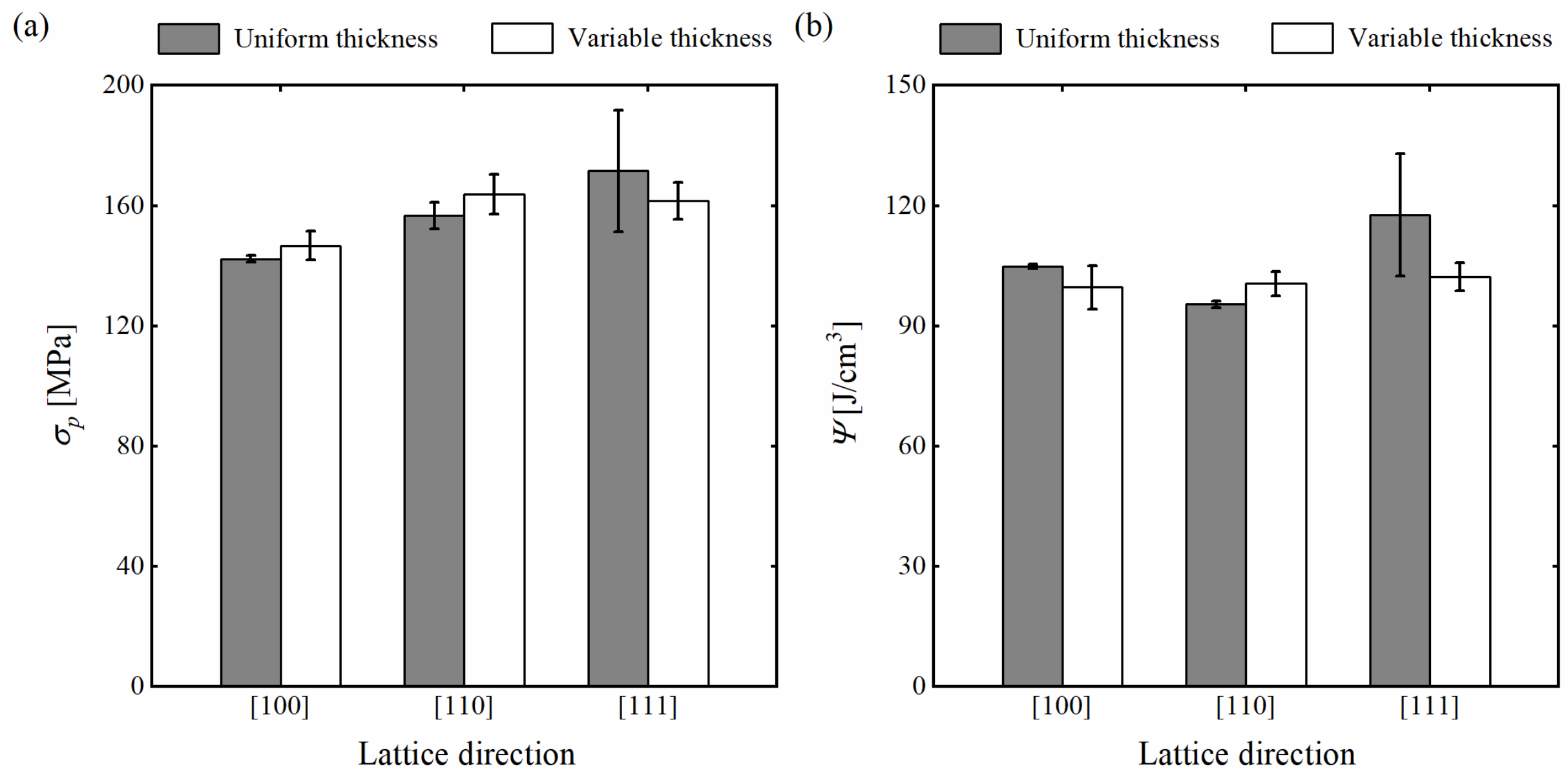

Fig. 16 (a) NPS, (b) SEA, along [100], [110] and [111] directions for the uniform and variable thickness specimens.

## 5. Conclusions

This work presents a family of elastically-isotropic open-cell variable thickness TPMS shell lattices with superior stiffness by a strain energy-based optimization algorithm. The cubic symmetry of the shell lattices is maintained by averaging the sensitivity terms in the optimization process. The optimization results show that all the six selected types of TPMS (N14, OCTO, IWP, FRD, N, P) shell lattices are able to achieve elastic isotropy by varying the shell thickness accordingly, among which N14 can maintain over 90% of the HS upper bound of bulk modulus. The Young's and shear moduli of the elastically-isotropic N14/OCTO/IWP/FRD/N lattices are found to reach nearly 40% ~ 60% of the HS upper bounds in low relative densities, which significantly outperform elastically-isotropic truss lattices of equal relative densities. Compared with the closed-cell elastically-isotropic plate lattices, the stiffness properties of the elastically-isotropic shell lattices are lower, as a trade-off induced by the open-cell topology. The maximum elastic properties achievable by open-cell lattices remain an open question for future research.

The uniform and variable thickness N14 shell lattices along [100], [110] and [111] directions are fabricated by micro-LPBF with SS316L and tested to verify the stiffness properties. The experimental results are in reasonable agreement with the FEA results. Comparisons of the tested Young's moduli along (1), (2), (6) and (7) directions reveal nearly elastic isotropy between [100] and [110] directions for the variable thickness specimens, with the moduli ratio reduced from 1.08 (uniform) to 1.02 (variable). Although the overall anisotropy in stiffness is still observed due to the anisotropy of the AM fabricated material, the elastic anisotropy of the lattices is shown to be reduced by varying the shell thickness accordingly. The following large strain compression tests reveal a layer-by-layer plastic buckling failure mode along [100] direction, while the failures along [110] and [111] directions are related to the shear deformation. It is also shown that the variable thickness specimens possess a reduced anisotropy of plateau stresses and nearly achieve isotropic energy absorption capacity.

In the experimental tests, the constituent material's anisotropy is observed, which induces anisotropy to the macroscopic stiffness properties of the variable thickness designs. A better way to take the anisotropy of the constituent material into consideration is to assume the constituent material to be orthotropic or transversely isotropic, instead of isotropic. Combination of this assumption with the whole analysis and optimization workflow in this work poses a potential opportunity to achieve a family of manufacturable elastically-isotropic open-cell shell lattices and will be explored in detail in our future work. Besides, the effects of thickness variations on the other properties of the lattices, such as permeability and fatigue performance, also deserve more explorations and will be studied through permeability measurements and fatigue tests in our future works.

## Declaration of Competing Interest

The authors declare that they have no known competing financial interests or personal relationships that could have appeared to influence the work reported in this paper.

## Acknowledgements

This work is supported by the Innovation and Technology Fund of the Government of the Hong Kong Special Administrative Region (Project No. ITS/008/19).

## Appendix A

The boundary conditions imposed on the 1/8 cell for the uniaxial strain loading ($\boldsymbol{\varepsilon}^{(U)} = [\varepsilon_{11}^{(U)}, 0, 0, 0, 0, 0]^T$) is:

$$\begin{aligned}
&(u_1^{x_0})^{(U)} = 0, \quad (\theta_2^{x_0})^{(U)} = 0, \quad (\theta_3^{x_0})^{(U)} = 0, \\
&(u_1^{x_+})^{(U)} = \varepsilon_{11}^{(U)} \frac{D}{2}, \quad (\theta_2^{x_+})^{(U)} = 0, \quad (\theta_3^{x_+})^{(U)} = 0, \\
&(u_2^{y_0})^{(U)} = 0, \quad (\theta_1^{y_0})^{(U)} = 0, \quad (\theta_3^{y_0})^{(U)} = 0, \\
&(u_2^{y_+})^{(U)} = 0, \quad (\theta_1^{y_+})^{(U)} = 0, \quad (\theta_3^{y_+})^{(U)} = 0, \\
&(u_3^{z_0})^{(U)} = 0, \quad (\theta_1^{z_0})^{(U)} = 0, \quad (\theta_2^{z_0})^{(U)} = 0, \\
&(u_3^{z_+})^{(U)} = 0, \quad (\theta_1^{z_+})^{(U)} = 0, \quad (\theta_2^{z_+})^{(U)} = 0,
\end{aligned} \tag{A1}$$

Similarly, those for the hydrostatic loading ($\boldsymbol{\varepsilon}^{(T)} = [\varepsilon_{11}^{(T)}, \varepsilon_{22}^{(T)}, \varepsilon_{33}^{(T)}, 0, 0, 0]^T$) is:

$$\begin{aligned}
&(u_1^{x_0})^{(T)} = 0, \quad (\theta_2^{x_0})^{(T)} = 0, \quad (\theta_3^{x_0})^{(T)} = 0, \\
&(u_1^{x_+})^{(T)} = \varepsilon_{11}^{(T)} \frac{D}{2}, \quad (\theta_2^{x_+})^{(T)} = 0, \quad (\theta_3^{x_+})^{(T)} = 0, \\
&(u_2^{y_0})^{(T)} = 0, \quad (\theta_1^{y_0})^{(T)} = 0, \quad (\theta_3^{y_0})^{(T)} = 0, \\
&(u_2^{y_+})^{(T)} = \varepsilon_{22}^{(T)} \frac{D}{2}, \quad (\theta_1^{y_+})^{(T)} = 0, \quad (\theta_3^{y_+})^{(T)} = 0, \\
&(u_3^{z_0})^{(T)} = 0, \quad (\theta_1^{z_0})^{(T)} = 0, \quad (\theta_2^{z_0})^{(T)} = 0, \\
&(u_3^{z_+})^{(T)} = \varepsilon_{33}^{(T)} \frac{D}{2}, \quad (\theta_1^{z_+})^{(T)} = 0, \quad (\theta_2^{z_+})^{(T)} = 0,
\end{aligned} \tag{A2}$$

while those for the pure shear loading ($\boldsymbol{\varepsilon}^{(S)} = [0, 0, 0, 0, 0, 2\varepsilon_{12}^{(S)}]^T$) is:

$$\begin{aligned}
&(u_2^{x_0})^{(S)} = 0, \quad (u_3^{x_0})^{(S)} = 0, \quad (\theta_1^{x_0})^{(S)} = 0, \\
&(u_2^{x_+})^{(S)} = \varepsilon_{12}^{(S)} \frac{D}{2}, \quad (u_3^{x_+})^{(S)} = 0, \quad (\theta_1^{x_+})^{(S)} = 0, \\
&(u_1^{y_0})^{(S)} = 0, \quad (u_3^{y_0})^{(S)} = 0, \quad (\theta_2^{y_0})^{(S)} = 0, \\
&(u_1^{y_+})^{(S)} = \varepsilon_{12}^{(S)} \frac{D}{2}, \quad (u_3^{y_+})^{(S)} = 0, \quad (\theta_2^{y_+})^{(S)} = 0, \\
&(u_3^{z_0})^{(S)} = 0, \quad (\theta_1^{z_0})^{(S)} = 0, \quad (\theta_2^{z_0})^{(S)} = 0, \\
&(u_3^{z_+})^{(S)} = 0, \quad (\theta_1^{z_+})^{(S)} = 0, \quad (\theta_2^{z_+})^{(S)} = 0,
\end{aligned} \tag{A3}$$

## Appendix B

For the sensitivity analysis of the optimization problem in Eq. (17), the Lagrangian function is first defined as:

$$L = J + \lambda g, \tag{B1}$$

where $\lambda$ denotes the Lagrange multiplier. The derivative of the objective function with respect to the shell thickness is derived as:

$$\frac{\partial J}{\partial \delta_i} = \frac{2(\xi - 1)}{(c_{11} - c_{12})^2}\left[\frac{\partial c_{44}}{\partial \delta_i}(c_{11} - c_{12}) - c_{44}\left(\frac{\partial c_{11}}{\partial \delta_i} - \frac{\partial c_{12}}{\partial \delta_i}\right)\right], \quad (i = 1, 2, ..., N), \tag{B2}$$

with the derivative terms $\frac{\partial c_{11}}{\partial \delta_i}$, $\frac{\partial c_{12}}{\partial \delta_i}$ and $\frac{\partial c_{44}}{\partial \delta_i}$ derived from Eq. (8a) ~ (8c) as:

$$\frac{\partial c_{11}}{\partial \delta_i} = \frac{2}{V}\frac{\partial e^{(U)}}{\partial \delta_i}, \quad (i = 1, 2, ..., N), \tag{B3a}$$

$$\frac{\partial c_{12}}{\partial \delta_i} = \frac{1}{V}\left[\frac{1}{3}\frac{\partial e^{(T)}}{\partial \delta_i} - \frac{\partial e^{(U)}}{\partial \delta_i}\right], \quad (i = 1, 2, ..., N), \tag{B3b}$$

$$\frac{\partial c_{44}}{\partial \delta_i} = \frac{2}{V}\frac{\partial e^{(S)}}{\partial \delta_i}, \quad (i = 1, 2, ..., N), \tag{B3c}$$

The derivatives of the strain energies in Eq. (B3a) ~ (B3c) are derived as [56]:

$$\frac{\partial e^{(j)}}{\partial \delta_i} = \frac{1}{2}[\boldsymbol{U}_i^{(j)}]^T \frac{\partial \boldsymbol{K}_i}{\partial \delta_i}\boldsymbol{U}_i^{(j)} = \frac{1}{2}\left\{[\boldsymbol{U}_{i(j)}^{(m)}]^T \frac{\partial \boldsymbol{K}_i^{(m)}}{\partial \delta_i}\boldsymbol{U}_{i(j)}^{(m)} + [\boldsymbol{U}_{i(j)}^{(b)}]^T \frac{\partial \boldsymbol{K}_i^{(b)}}{\partial \delta_i}\boldsymbol{U}_{i(j)}^{(b)}\right\}, \quad (i = 1, 2, ..., N; \quad j = U, T, S), \tag{B4}$$

where $\boldsymbol{U}_i$ and $\boldsymbol{K}_i$ denote the element displacement vector and element stiffness matrix, respectively. Eq. (B4) reveals that $\frac{\partial e^{(j)}}{\partial \delta_i}$ can be decomposed into two parts, i.e., the membrane term and bending term, in which $\boldsymbol{U}_i^{(m)}$ and $\boldsymbol{U}_i^{(b)}$ denote the re-organized element membrane displacement vector and element bending displacement vector, while $\boldsymbol{K}_i^{(m)}$ and $\boldsymbol{K}_i^{(b)}$ denote the re-organized element membrane stiffness matrix and element bending stiffness matrix. The element membrane stiffness matrix is derived as:

$$\boldsymbol{K}_i^{(m)} = \int_{M_i} [\boldsymbol{B}^{(m)}]^T \boldsymbol{D}^{(m)} \boldsymbol{B}^{(m)} dA, \quad (i = 1, 2, ..., N),$$

$$\boldsymbol{D}^{(m)} = \frac{E_s \delta_i}{1 - \nu_s^2} \begin{pmatrix} 1 & \nu_s & 0 \\ \nu_s & 1 & 0 \\ 0 & 0 & \frac{1 - \nu_s}{2} \end{pmatrix}, \tag{B5}$$

where $\boldsymbol{B}^{(m)}$ denotes the membrane strain-displacement matrix, while $\delta_i$ and $M_i$ denote the thickness and mid-surface of element $i$. The element bending stiffness matrix is derived from the Kirchhoff-Love plate/shell theory as [57, 58]:

$$\boldsymbol{K}_i^{(b)} = \int_{M_i} [\boldsymbol{B}^{(b)}]^T \boldsymbol{D}^{(b)} \boldsymbol{B}^{(b)} dA, \quad (i = 1, 2, ..., N),$$

$$\boldsymbol{D}^{(b)} = \frac{E_s \delta_i^3}{12(1 - \nu_s^2)} \begin{pmatrix} 1 & \nu_s & 0 \\ \nu_s & 1 & 0 \\ 0 & 0 & \frac{1 - \nu_s}{2} \end{pmatrix}, \tag{B6}$$

where $\boldsymbol{B}^{(b)}$ denotes the bending strain-displacement matrix. Since both the strain-displacement matrices $\boldsymbol{B}^{(m)}$ and $\boldsymbol{B}^{(b)}$ are independent of the shell thickness $\delta_i$, the derivatives of the membrane and bending stiffness matrices in Eq. (B4) are thus obtained as:

$$\frac{\partial \boldsymbol{K}_i^{(m)}}{\partial \delta_i} = \frac{1}{\delta_i} \boldsymbol{K}_i^{(m)}, \quad (i = 1, 2, ..., N), \tag{B7a}$$

$$\frac{\partial \boldsymbol{K}_i^{(b)}}{\partial \delta_i} = \frac{3}{\delta_i} \boldsymbol{K}_i^{(b)}, \quad (i = 1, 2, ..., N), \tag{B7b}$$

The element membrane and bending stiffness matrices $\boldsymbol{K}_i^{(m)}$ and $\boldsymbol{K}_i^{(b)}$ are extracted from ABAQUS by adding commands in the .inp file, based on which the derivative terms in Eq. (B4) and (B7a) ~ (B7b) are evaluated. The derivatives of the objective function are then achieved by Eq. (B2) and (B3a) ~ (B3c). On the other hand, the derivative of the constraint function with respect to the shell thickness is derived as:

$$\frac{\partial g}{\partial \delta_i} = \frac{1}{V} m_i, \quad (i = 1, 2, ..., N), \tag{B8}$$

where $m_i$ denotes the mid-surface area of element $i$, based on which the derivative of the Lagrangian function is derived as:

$$\frac{\partial L}{\partial \delta_i} = \frac{\partial J}{\partial \delta_i} + \lambda \frac{\partial g}{\partial \delta_i}, \quad (i = 1, 2, ..., N), \tag{B9}$$

To obtain the updating scheme of the shell thickness, the derivative of the Lagrangian function with respect to the pseudo-time is derived by the chain rule as:

$$\dot{L} = \sum_{i=1}^{N} \frac{\partial L}{\partial \delta_i} \dot{\delta}_i, \tag{B10}$$

In this work, the steepest descent method is adopted as the optimizer, meaning $\dot{\delta}_i$ is taken as:

$$\dot{\delta}_i = -\frac{\partial L}{\partial \delta_i}, \quad (i = 1, 2, ..., N), \tag{B11}$$

so that $\dot{L} = -\sum_{i=1}^{N} (\frac{\partial L}{\partial \delta_i})^2 \leq 0$. The element thickness $\delta_i$ is then updated by:

$$\delta_i(a+1) = \delta_i(a) + dt \cdot \dot{\delta}_i(a), \quad (i = 1, 2, ..., N; a \geq 1), \tag{B12}$$

where $a$ and $dt$ denote the iteration step and step size, respectively. In each iteration step, the Lagrange multiplier is updated by the augmented Lagrangian method as [59]:

$$\begin{aligned} &\lambda_{a+1} = \lambda_a + \mu_a g_a, \quad (a \geq 1), \\ &\mu_{a+1} = \min(\mu_a + \Delta\mu, \mu_{\max}), \quad (a \geq 1), \end{aligned} \tag{B13}$$

where $\Delta\mu$ and $\mu_{\max}$ denote the increment and maximum value of the parameter $\mu$, respectively.